\documentclass{amsart}
\usepackage{amsmath} 
\usepackage{amssymb}
\usepackage{mathrsfs}

\newtheorem{theorem}{Theorem}[section] 
\newtheorem{claim}[theorem]{Claim}

\theoremstyle{definition}
\newtheorem{definition}[theorem]{Definition}

\newtheorem{discussion}[theorem]{Discussion}

\newtheorem{question}[theorem]{Question}

\theoremstyle{remark}
\newtheorem{remark}[theorem]{Remark}
\newtheorem{notation}[theorem]{Notation}
\newtheorem{conclusion}[theorem]{Conclusion}

\newcommand{\rest}{{\restriction}}
 
\newcommand{\Dom}{{\rm Dom}} 

\newcommand{\otp}{{\rm otp}}
\newcommand{\nor}{{\rm nor}}

\newcommand{\id}{{\rm id}}
\newcommand{\Ln}{{\rm Ln}}

\newcommand{\stg}{{\rm stg}}

\newcommand{\Suc}{{\rm Suc}}

\newcommand{\Sub}{{\rm Sub}}
\newcommand{\Lim}{{\rm Lim}}
\newcommand{\Fun}{{\rm Fun}}
\newcommand{\Ord}{{\rm Ord}}

\newcommand{\unf}{{\rm unf}}
\newcommand{\cf}{{\rm cf}}

\newcommand{\stat}{{\rm stat}}
\newcommand{\wilog}{{\rm without loss of generality}}

\newcommand{\then}{{\underline{then}}}
\newcommand{\when}{{\underline{when}}}
\newcommand{\Then}{{\underline{Then}}}

\newcommand{\Iff}{{\underline{iff}}}
\newcommand{\mn}{{\medskip\noindent}}
\newcommand{\sn}{{\smallskip\noindent}}

\newcommand{\gK}{{\mathfrak K}}

\newcommand{\gf}{{\mathfrak f}}
\newcommand{\gq}{{\mathfrak q}}

\newcommand{\cJ}{{\mathcal J}}

\newcommand{\bbD}{{\mathbb D}}
\newcommand{\cD}{{\mathcal D}}

\newcommand{\cH}{{\mathcal H}}

\newcommand{\cF}{{\mathcal F}}

\newcommand{\cI}{{\mathcal I}}

\newcommand{\bbZ}{{\mathbb Z}}

\newcommand{\bbP}{{\mathbb P}}
\newcommand{\cP}{{\mathcal P}}
\newcommand{\gp}{{\mathfrak p}}

\newcommand{\bbQ}{{\mathbb Q}}

\newcommand{\cS}{{\mathcal S}}
\newcommand{\cT}{{\mathcal T}}
\newcommand{\cW}{{\mathcal W}}
\newcommand{\gt}{{\mathfrak t}} 
\newcommand{\gT}{{\mathfrak T}} 
\newcommand{\gs}{{\mathfrak s}}

\newcommand{\cY}{{\mathcal Y}}

\newcount\skewfactor
\def\mathunderaccent#1#2 {\let\theaccent#1\skewfactor#2
\mathpalette\putaccentunder}
\def\putaccentunder#1#2{\oalign{$#1#2$\crcr\hidewidth
\vbox to.2ex{\hbox{$#1\skew\skewfactor\theaccent{}$}\vss}\hidewidth}}
\def\name{\mathunderaccent\tilde-3 }

\newenvironment{PROOF}[2][\proofname.]
   {\begin{proof}[#1]\renewcommand{\qedsymbol}{\textsquare$_{\rm #2}$}}
   {\end{proof}}

\begin{document}

\title {Many forcing axioms for all regular
  uncountable cardinals}
\author {Saharon Shelah}
\address{Einstein Institute of Mathematics\\
Edmond J. Safra Campus, Givat Ram\\
The Hebrew University of Jerusalem\\
Jerusalem, 91904, Israel\\
 and \\
 Department of Mathematics\\
 Hill Center - Busch Campus \\ 
 Rutgers, The State University of New Jersey \\
 110 Frelinghuysen Road \\
 Piscataway, NJ 08854-8019 USA}
\email{shelah@math.huji.ac.il}
\urladdr{http://shelah.logic.at}
\thanks{The author thanks Alice Leonhardt for the beautiful typing.
Research supported by the German-Israeli Foundation for Scientific
Research and Development (Grant No:I-706054.6/2001)  First typed June
2, 2002. Publication 832}






\date{June 21, 2013}

\begin{abstract}
Our original aim was, in Abelian group theory to 
prove the consistency of: $\lambda$ is strong limit
singular and for some properties of abelian groups which are relatives
of being free, the compactness in singular fails.  In fact this should
work for $R$-modules, etc.  As in earlier cases part of the work is analyzing
how to move between the set theory and the algebra.   

Set theoretically we try to force a universe which satisfies G.C.H. and diamond
holds for many stationary sets but, for every regular uncountable
$\lambda$, in some sense 
anything which ``may" hold for some stationary set, does
hold for some stationary set.  More specifically we try to get a universe
satisfying GCH such that e.g. for regular $\kappa < \lambda$ there are
pairs $(S,B),S \subseteq S^\lambda_\kappa$ stationary, $B \subseteq
{\cH}(\lambda)$, which satisfies some pregiven forcing axiom
related to $(S,B)$, (so $(\lambda \backslash S)$-complete,
i.e. ``trivial outside $S$) \underline{but} no more, i.e. slightly stronger
versions fail.  So set theoretically we try to get a universe
satisfying G.C.H. but still satisfies ``many", even for a maximal
family in some sense, of forcing axioms of the form ``for some
stationary" while preserving GCH.  As completion of the work lag for
long, here we deal only with the set theory.
\end{abstract}

\maketitle
\numberwithin{equation}{section}
\setcounter{section}{-1}
\newpage

\centerline{Annotated Content} 
\bigskip

\noindent
\S0 \quad Introduction, pg.\pageref{Introduction}
\bigskip

\noindent
\S1 \quad The iteration on all cardinals, pg.\pageref{Theit}
\mn
\begin{enumerate}
\item[${{}}$]   [We define when ${\gt}$ is a $\lambda$-task
template, $(< \alpha)$-strategically
$S$-complete, $\lambda$-tasks and cases of $\lambda$-task templates
and when a $\lambda$-task is satisfied, i.e. each instance is (see
\ref{1t.1}, \ref{1t.5}).  
We prove that if $\lambda = \lambda^{< \lambda} > \aleph_0$ and 
$2^\lambda = \lambda^+$ then there is a suitable forcing $\bbP$ 
taking care of all $\lambda$-tasks (\ref{1t.13}) and that we can
do it for all cardinals, (\ref{1t.15}).]
\end{enumerate}
\bigskip

\noindent
\S2 \quad An example: Relatives of diamonds, pg.\pageref{anex}
\bigskip

\noindent
\S3 \quad Parameters for completeness of forcing, pg.\pageref{Parameters}
\mn
\begin{enumerate}
\item[${{}}$]   [Also this section is set theoretic but fine tuned to
  uniformization and hence Abelian groups (or modules) problems.
We define when $(D,{\cT})$ is a $(\mu,\kappa,\theta)$-special pair, 
(\ref{2f.3}) and when $\bold p = (D,{\cT},\lambda,\mathscr{S},\bar \eta)$ is a
$(\mu,\kappa,\theta)$-parameter (\ref{2f.7}(1)) and prove iteration claims for
$\bold p$ (in \ref{2f.11}).  We then look for suitable filters (see
Definition \ref{2f.13}), and prove existence (\ref{2f.15},\ref{2f.17}).
We point out the case of stationary $\cS \subseteq {}^{\delta
\ge}([\alpha]^{< \mu})$, $\delta$-fat $S \subseteq
\lambda$ (i.e., $\cS = {}^{\delta \ge} S$) and $\cS$
complete, though actually our main forcing are
strategically complete with the $\bold p$ being forced.
With those notions
we may define uniformization properties and $\bold p$-completeness of
forcing notions (\ref{2f.9}(2)).  This will allow us to iterate
\ref{2f.11}(1).   We note that there is
$(\mu,\mu,\aleph_0)$-special pairs if $\diamondsuit_\mu$
(\ref{2f.15}), and that we can combine
$(\mu_\delta,\kappa,\theta)$-special pairs for enough $\delta$'s to a
$(\mu,\kappa,\theta)$-special pair (see \ref{2f.17}).]
\end{enumerate}
\newpage

\section {Introduction} \label{Introduction}

The analysis of those questions lead to having many uniformation
properties to many regular cardinals smaller then the singulars.
\bigskip

\subsection  {An Abelian group theory motivation} \label{Anabelian} \
\bigskip

Compactness in singular (provable in ZFC) play a crucial role in
the solution of Whitehead problem, see \cite{Sh:52},\cite{EM02}.

There are some ``compactness in singular cardinals" theorems for
Abelian groups when we assume
$\bold V = \bold L$ (or something in this direction), of Eklof and more
recently of Struengman.  A natural
question is whether the $\bold V = \bold L$ is needed.  In fact the
proofs of each of those theorems can be decomposed to 
two parts: one part which is
compactness in singular for another condition (say being a free
Abelian group); another part is the equivalence of the two conditions 
using $\bold V = \bold L$.  Anyhow this stands behind the following questions:

\begin{question}
\label{z3} 
Compactness in singulars 
for $\{G:\text{Ext}(G,\bbZ) = 0\}$, 
which arise from an old work of Eklof, \cite[Theorem 8.5]{Ek80}, 
see more in \cite{EFSh:352}.
\end{question}

\begin{question}
\label{z6}
Does compactness in singulars for
$\{G:\text{Ext}(G,T)=0\}$ holds (for $T$ any torsion group)?

The question above, \ref{z6}, 
was asked by Struengmann following the paper \cite[Proposition
2.6]{Str02} and was the immediate reason for this work on \ref{z3},
\ref{c6} and more.  Later Eklof has asked me on \ref{z6}.  We shall show the
consistency of the negation in both cases.
\end{question}

\begin{notation}
\label{z9}
For $C$ a set of ordinals, 
\mn
\begin{enumerate}
\item[$(a)$]   let $(\forall^* \alpha \in C)\varphi$ means $\{\alpha \in
C:\varphi(\alpha)\}$ contains a co-bounded subset of $C$ 
\sn
\item[$(b)$]   let $(\forall^D \alpha \in C)\varphi(\alpha)$ 
means $\{\alpha \in C:\varphi(\alpha)\} \in D$, for $D$ a filter on $C$
\sn
\item[$(c)$]   using $\forall^I,I$ an ideal on $C$ means
$\forall^D,D$ the dual filter.
\end{enumerate}
\end{notation}
\bigskip

\subsection {The set theoretic view} \label{Theset} \
\bigskip

As said above our original motivation concerns Abelian groups and
consistency of incompactness in singulars for algebraic problem we
explained above.  Explain the set theoretic side, starting from
uniformization (see \cite{Sh:587} and history there).  Consider
$\lambda$ regular uncountable, stationary $S \subseteq \lambda$
consisting of limit ordinals; ladder system $\bar C = \langle
C_\delta:\delta \in S \rangle$, i.e., $C_\delta \subseteq \delta =
\sup(C_\delta)$ (with $C_\delta$ not necesarily a closed subset of
$\delta$) and $\bold h:\lambda
\rightarrow \lambda$, (if $\bold h$ is constantly 2 we may omit it; we
do not consider here the case $f_\delta \in {}^{(C_\delta)}
\delta$)).  Here and later we may replace $(\forall \delta \in S)$ by
$(\forall^{\cD_\lambda} \delta \in S)$ where $\cD_\lambda$ is e.g. the
club filter on $S$, see \ref{z9}; 
presently this does not cause a great difference,
but the property is weaker.
\mn
\begin{enumerate}
\item[$\circledast_1$]  we say that $\bar C$ has $\bold
h$-uniformization \when \, for every sequence $\langle f_\delta:\delta \in
S\rangle$ satisfying $f_\delta \in \prod\limits_{\alpha \in C_\delta} 
\bold h(\delta)$ for $\delta \in S$ there is $f \in
{}^\lambda \lambda$ such that $(\forall \delta \in S)(\forall^* \alpha
\in C_\delta)(f(\alpha)) = f_\delta(\alpha)$. 
\end{enumerate}
\mn
We may add $\bar D = \langle D_\delta:\delta \in S \rangle,D_\delta$ a
filter on $C_\delta$ and consider
\mn
\begin{enumerate}
\item[$\circledast_2$]   we say that $\bar C$ has $(\bold h,\bar
D)$-uniformization if for every $\bar f = \langle f_\delta:\delta \in
S \rangle,f_\delta = \prod\limits_{\alpha \in C_\alpha} \bold
h(\alpha)$ there is $f \in \prod\limits_{\alpha < \lambda} \bold
h(\alpha)$ such that $(\forall \delta \in S)(\forall^D \alpha \in
C_\delta)f(\alpha) = f_\delta(\alpha)$, i.e. $\{\alpha \in
C_\delta:f(\alpha) = f_\delta(\alpha)\} \in D_\delta$ for $\delta \in S$.
\end{enumerate}
\mn
We use $S(\bar C) = S$.

We may consider questions close to Abelian group theory
\mn
\begin{enumerate}
\item[$\circledast_3$]   (group-uniformity) for a sequence $\bar K = 
\langle K_\alpha:\alpha < \lambda\rangle$ of groups and $\bar C$ as above let
$K^\delta = \prod\limits_{\alpha \in C_\delta} K_\alpha$ 
and let ps-ext$(\bar K,\bar C,D)$ be
$K^* / \unf(\bar K,\bar C,\bar D)$ where $K_* =
\prod\limits_{\delta \in S} K^\delta$ and $\unf(\bar K,\bar C,\bar D) 
= \{f \in S(\bar C)$: there is $h \in 
\prod\limits_{\alpha < \lambda} K_\alpha$ 
such that $(\forall \delta \in S)(\forall^\infty
\alpha \in C_\delta)[f_\delta(\alpha) = f(\alpha)]$, i.e. $\{\alpha \in
C_\delta:f(\alpha) = f_\delta(\alpha)\} \in D_\delta$ for $\delta \in S\}$.
\end{enumerate}
\mn
We may vary more: for some function $\bold F$ replace $f(\alpha) =
f_\delta(\alpha)$ by $f_\delta(\alpha) = \bold F(f \restriction
C_\delta \cap (\alpha +1)$, this is closer to Ext.

We concentrate on the case that G.C.H. holds and we have uniformation
for enough stationary sets $S$ for some appropriate $\bar C$.  But
here we try to do it for every regular uncountable $\lambda$ and for
all ``tasks" of such forms.  On earlier works forcing uniformization
see Eklof-Mekler \cite{EM02}.

So usually we have to assume:
\mn
\begin{enumerate} 
\item[$\circledast_4$]   $\lambda = \mu^+$ and $\delta \in S
\Rightarrow \text{ cf}(\delta) = \text{ cf}(\mu)$.
\end{enumerate}
\mn
We can force (using relatives of pseudo-completeness, see \S2)
\mn
\begin{enumerate} 
\item[$\circledast_5$]  $(a) \quad$ if a stationary $S \subseteq \lambda$
satisfies $\circledast_4$ then for some stationary

\hskip25pt  $S' \subseteq S$
some ladder system $\langle C_\delta:\delta \in S' \rangle$
has uniformization
\sn
\item[${{}}$]  $(b) \quad$ even if $\bar C = \langle C_\delta:\delta
  \in S\rangle$ is a ladder system then for some stationary

\hskip25pt  $S' \subseteq S,\bar C \rest S$ has uniformization.
\end{enumerate}
\mn
Here we shall be interested in getting distinction between quite close
relatives of this, so we have to force ``less" and the problem is to
phrase exactly what we like to have but not
to get more uniformization than we intend.  In \cite{Sh:587} we
consider, for other reasons the case $\lambda = \mu^+,\mu = 
\cf(\mu),D_\delta$ a somewhat regular filter on $C_\delta$ (so
$|C_\delta| = \mu$).   The aim there was to get some uniformization on
$S^\lambda_\mu$, this could have been done also 
in \cite{Sh:667} which concentrates on successor of singulars but there was no
point gain by it.  The iteration in \S1 preserves G.C.H. and for
every successor $\lambda = \text{ cf}(\lambda) > \aleph_0$ and stationary $S
\subseteq \lambda$ (mainly satisfying $\circledast_4$), add a 
stationary subset $S' \subseteq S$ and on it for
something which is similar enough to a case of uniformization; so called
$\lambda$-task, which may be relevant for 
Abelian groups.  This is arranged such that there is little
interaction between forcing for the 
different tasks.  For each regular uncountable
$\lambda$ there is an iterated forcing which adds enough ``tasks" and mainly
add solutions for each case of each such task.  There is no problem in
the forcing because if the task is ``too hard" (e.g. gives a provably
impossible situation), the forcing will, e.g. make
$S'$ not be stationary (for abelian group - the group we add becomes
free or examples intended to show Ext$(G,\Bbb Z) \ne 0$ are no longer so).  So
as in many cases ``we solve our problems by putting them on someone
else's shoulders".

The central case here is with $D_\delta$ a regular filter; then some
relevant finitary linear combinations become critical.
A major part
in the proof is trying to prove that the iterated forcing gives ``good" limit.
As we arrange it, not collapsing cardinals holds trivially (as we just
ask ``for some stationary $S' \subseteq S$").  So the difficulty is 
preserving stationarity of relevant sets, and/or showing that 
undesirable objects are not added by building inverse system of trees
of forcing conditions, (a central method in such consistency proofs).  In the 
main case, inside the forcing proof we have to
consider only finitely many coordinates, using the regularity of the
relevant filter.  So trying to immitate the
argument in \cite{Sh:125}, we can do it in higher cardinals, if we
carry with us strong enough induction hypothesis.  Should we in the
forcing in \S1 add non-reflecting stationary sets?  We may use such
$\lambda$-tasks, but we may like to preserve supercompactness, so we
do not like to, but then we have to use reflection.

As finishing the full work large for too long, we delay the results on Ent.
This work and the application to Ent were 
presented at the meeting in honor of Eklof in Summer 2008.
\bigskip

\subsection {Preliminaries} \label{Prelim} \
\bigskip

\begin{definition}
\label{z21}
1) For $\lambda$ regular uncountable and $\cS \subseteq \cP(\lambda)$
   let $\nor-\id(\cS)$ be the minimal normal ideal which includes $\cS$.
\end{definition}

\begin{definition}
\label{z23}
Let $\bbD$ be a filter on ${}^{\lambda >}(\cP(\lambda))$.

\noindent
1) We call $\bbD$ normal \when \, for every $\chi > \lambda$ and $x
\in \cH(\chi)$ there is $\cY \in \bbD$ such that: $\langle N_i \cap
\lambda:i \le \delta\rangle \in \cY$ wherever $\bar N = \langle N_i:i
\le \delta\rangle$ obeys $(\bbD,\chi,\lambda,x)$ 
which means (omitting $\cY$ means for some $\cY$)
\mn
\begin{enumerate}
\item[$(a)$]  $N_i \prec (\cH(\chi),\in)$
\sn
\item[$(b)$]  $\lambda,\bbD \in N_i$ 
\sn
\item[$(c)$]  $x \in N_i$  
\sn
\item[$(d)$]  $N_i$ is increasing continuous
\sn
\item[$(e)$]  $\bar N \rest (i+1) \in N_{i+1}$.
\end{enumerate}
\mn
2) We say $\bbD$ is a fat normal (filter) \when \, for every
$\alpha < \lambda$ we have $\cY \in \bbD \wedge \alpha < \lambda
\Rightarrow \cY_{\ge \alpha} := \{\bar u \in \cY:\ell g(\bar u) >
\alpha\} \ne \emptyset \mod \bbD$.

\noindent
3) We say a forcing notion $\bbP$ is $\bbD$-complete \when \,:
\mn
\begin{enumerate}
\item[$(a)$]  forcing by $\bbP$ adds no new sequence of ordinals of
  lenth $< \lambda$
\sn
\item[$(b)$]  for every $\chi$ for some $Y \in \bbD$ and 
$x \in \cH(\chi)$, for every $\bar N$ obeying $(\bbD,\chi,\lambda,x)$ we have:
\sn
\begin{enumerate}
\item[$\bullet$]  if $j \le \ell g(\bar N)$ is a limit ordinals, 
$\bar p = \langle p_i:i < j \rangle \in \prod\limits_{i < j} 
(N_i \cap \bbP)$ is increasing $i<j \Rightarrow p \rest (i+1) \in N_i$
and weakly generic for $\bar N$ (if $i<j$ then for some $n$, for every
$\bbP$-name of an ordinal $\name\tau \in N_i,p_{i+n} \Vdash
``\name\tau \in N_{i+n} \cap \Ord"$, as in \cite{Sh:587}) \then \,
$\bar p$ has a $\le_{\bbP}$-upper bound (hence preserve 
${}^{\lambda>} \bold V$).
\end{enumerate}
\end{enumerate}
\mn
4) For stationary $w \subseteq \lambda = \text{ cf}(\lambda)$ let
$\bbD_{\lambda,w}$ be the minimal normal filter on ${}^{\lambda
  >}\cP(\lambda)$ to which the following set belongs
$\{\bar\alpha:\bar\alpha$ is an increasing continuous sequence 
of ordinals from $w\}$.

\noindent
5) If $\bbD$ is normal, let prog$_{\bbD}(\chi) = \{\cY:\cY \subseteq
   \{N:N \prec \cH(\chi),\in)$ and $\|N\| < \lambda\}$ and for some $x \in
\cH(\chi)$: the set $\cY$ includes $\{N_\delta$: there is $\bar N =
\langle n_i:i \le \delta \rangle$ which obeys $(\bbD,\chi,\lambda,x)\}\}\}$.
\end{definition}

\begin{definition}
\label{z28}
Let $\Fun(G,H)$ be the set of functions from $G$ to $H$, not
necessarily homomorphisms.
\end{definition}
\bigskip

\centerline {$* \qquad * \qquad *$}
\bigskip

Concerning strategic completeness
\begin{definition}
\label{a50}
1) For a regular uncountable $\lambda$ and ordinal $\xi \le \lambda$
   we say a forcing notion $\bbQ$ is 
$(\xi,\bold S_1,\name{\cS}_2)-\stg$-complete \when \,:
\mn
\begin{enumerate}
\item[$(a)$]  $\cS_2$ is a family of subsets of $\lambda$
\sn
\item[$(b)$]  $\name{\cS}_2$ is a $\bbQ$-name of a family of subsets
  of $\lambda$
\sn
\item[$(c)$]  in the game $\Game =
  \Game_{\alpha_*}(\cS_1,\name{\cS}_2,\bbQ)$ the completeness player
  has a winning strategy
\sn
\item[${{}}$]  $(\alpha) \quad$ a play last $\alpha_*$ moves
\sn
\item[${{}}$]  $(\beta) \quad$ in the $\alpha$-th move a quadruple
  $(p_\alpha,S_{1,\alpha},\name S_{2,\alpha},\varepsilon_\alpha)$ is
  chosen such that
\sn
\begin{enumerate}
\item[${{}}$]  $\bullet_1 \quad p_\alpha \in \bbQ$ such that $\beta <
  \alpha \Rightarrow p_\beta \le_{\bbQ} p_\alpha$
\sn
\item[${{}}$]  $\bullet_2 \quad \varepsilon_\alpha < \lambda$ 
such that $\langle \varepsilon_\beta:\beta \le \alpha\rangle$ is
increasing continuous
\sn
\item[${{}}$]  $\bullet_3 \quad S_{1,\alpha} \in
  \nor-\id_\lambda(\cS_1)$ such that $\langle S_{1,\beta}:\beta \le 
\alpha \rangle$ is $\subseteq$-increasing 

\hskip30pt continuous
\sn
\item[${{}}$]  $\bullet_4 \quad \name\cS_{2,\alpha} \in
  \nor-\id_\lambda(\name{\cS}_2)$ such that $\langle \name
  S_{2,\beta}:\beta \le \alpha \rangle$ is $\subseteq$-increasing

\hskip30pt   continuous
\sn
\item[${{}}$]  $\bullet_5 \quad p_\alpha \Vdash ``\varepsilon_\alpha
  \notin S_{1,\alpha} \cup \name S_{2,\alpha}$
\end{enumerate}
\sn
\item[${{}}$]  $\gamma) \quad$ if $\alpha = 2n$ or $\alpha = \omega(1
  + \beta) + 2n +1$ \then \, the incompleteness player

\hskip25pt  chooses the
  quadruple, otherwise the completeness player chooses
\sn
\item[${{}}$]  $(\delta) \quad$ the completeness player wins a play
  \Iff \, he has a legal move for every 

\hskip25pt  limit $\alpha < \alpha_*$ (the only problematic point is choosing).
\end{enumerate}
\end{definition}

\noindent
We have natural iteration claims.
\begin{claim}
\label{a52}
If (A) then (B) where:
\mn
\begin{enumerate}
\item[$(A)$]  $\gq$ is a $(\lambda,\xi)-\stg$-iteration which means
\sn
\begin{enumerate}
\item[$(a)$]  $\lambda = \lambda^{< \lambda} > \aleph_0$ and $\xi \le
  \lambda$
\sn
\item[$(b)$]  $\gq = \langle
  \bbP_\alpha,\bbQ_\beta,\name{\cS}_\beta:\alpha \le \alpha_*,\beta <
  \alpha_*$
\sn
\item[$(c)$]  $\langle \bbP_\alpha,\name{\bbQ}_\beta:\alpha \le
  \alpha_*,\beta < \alpha_*\rangle$ is $(\le \lambda)$-support
  iteration with limit $\bbP_{\alpha_*}$
\sn
\item[$(d)$]  $\name{\cS}_\alpha$ is a $\bbP_\alpha$-name of a
  $\subseteq$-increasing family of subsets of $\lambda$
\sn
\item[$(e)$]  $\Vdash_{\bbP_\alpha} ``\lambda \notin
  \nor-\id_\lambda(\name{\cS}_\alpha)"$
\sn
\item[$(f)$]  $\Vdash_{\bbP_\alpha} ``\bbQ_\alpha$ is
  $(\xi,\name{\cS}_\alpha,\name{\cS}_{\alpha+1})-\stg$-complete
\end{enumerate}
\sn
\item[$(B)$]  $(a) \quad \bbP_{\alpha_*}$ is
  $(\xi,\cS_0,\name{\cS}_{\alpha_*})-\stg$-complete
\sn
\item[${{}}$]  $(b) \quad$ even $\bbP_\beta/\bbP_\alpha$ is when
$\alpha < \beta \le \alpha_*$.
\end{enumerate} 
\end{claim}

\begin{PROOF}{\ref{a52}}
Straight.
\end{PROOF}

\begin{claim}
\label{a55}
We have $\gq$ is a $(\lambda,\xi)-\stg$-iteration \when \,:
\mn
\begin{enumerate}
\item[$(A)'$]  clauses (a),(c),(d) of \ref{a52}(A)
\sn
\item[${{}}$]  $(\alpha) \quad \lambda \notin \nor-\id_\lambda(\cS_0)$
\sn
\item[${{}}$]  $(\beta) \quad$ if $\alpha < \alpha_{gq}$ and
  $\Vdash_{\bbP_\alpha}$ ``if $\lambda \notin
  \nor-\id_\lambda(\name{\cS}_\alpha)$ then the forcing notion

\hskip25pt $\name{\bbQ}_\alpha$ is $(\xi,\name{\cS}_\alpha,
\name{\cS}_{\alpha +1})-\stg$-complete
\sn
\item[${{}}$]  $(\gamma) \quad$ if $\alpha$ is a limit ordinal $\le
  \alpha_{\gq}$ then $\name{\cS}_\alpha =
  \cup[\{\name{\cS}_\beta:\beta < \alpha\}$ or 

\hskip25pt  $\lambda \notin \nor-\id_\lambda(\cS_\beta)$.
\end{enumerate}
\end{claim}

\begin{remark}
We may weaken $(\gamma)$ in $(A)'$ of \ref{a55} , but no use for now.
\end{remark}
\newpage

\section {The iteration on all cardinals} \label{Theit}

This section is purely set theoretical, this is continued in \S2 but
there we deal with iteration, whereas there it deals with ``one forcing".

Our program is as follows, we start with $\bold V$ satisfying
G.C.H. and preserve it.  By the character of our aims we have to have
relevant existence in all cardinals.
The freedom we have is to try for each regular uncountable
$\mu$ to add a stationary subset $S$ of $\mu$ for each task, and try
to fulfill the task.  Now if the task is ``too hard", say it
contradicted by ZFC + GCH, the iteration may make $S$ not stationary
(or ``harm" it in other ways).  So in our frame this causes no problem
here, it just passes the burden to proofs about the specific tasks.

We may in the forcing axioms (related to specific stationary $S
\subseteq M$) replace ``of cardinality $\mu$" by ``of
cardinality $\mu^+$" satisfying a strong form of $\mu^+$-c.c. but this
is not central here.
The forcing for Abelian groups in $\lambda = \mu^+$ require we have
such abelian groups in $\mu^+$ and only stationary $S \subseteq
\{\delta < \lambda:\cf(\delta) = \cf(\mu)\}$ are used, there are
differences according to $\lambda$ being successor of regular,
successor of singular or inaccessible.

Now \ref{1t.1} defines a $\lambda$-task template.  We intend for
each $\lambda$-task to force an example, then later in the iteration we
force a solution to each relevant case; 
we hope the rest of the iteration do preserve some
desired properties; this means that for any $B \subseteq 
{\cH}(\lambda)$ we have to deal with each specific
$\lambda$-tasks derived from it (see Definition \ref{1t.5}).

\begin{definition}
\label{1t.1} 
For $\lambda$ regular uncountable
we say ${\gt}$ is a $\lambda$-task template \Iff \, ${\gt}$ consists of:
\mn
\begin{enumerate}
\item[$(a)$]   $S^{\gt}$, a subset of $\lambda$, normally stationary 
\sn
\item[$(b)$]    $\bbQ_{\gt} = \bbQ^{\gt} = \bbQ[{\gt}]$, a forcing
notion of cardinality $\lambda$ which is $(< \lambda)$-strategically
complete, see \ref{1t.3}; [not a serious difference if we ask 
$\lambda$-strategically complete but see (c)$(\beta)$]
\sn
\item[$(c)$]   a function $p \mapsto S_p = S^{\gt}_p$ 
from $\bbQ_{\gt}$ to $\{S:S$ a bounded subset of $S^{\gt}$
presented as a characteristic function from some ordinal $< \lambda$ to
$\{0,1\}\}$ and let $\delta_{\gt}(p) := \Dom(S^{\gt}_p)$ so an
ordinal $< \lambda$, such that:
\sn
\begin{enumerate}
\item[$(\alpha)$]   $\bbQ_{\gt} \models p \le q \Rightarrow S_p = 
S_q \cap \delta_{\gt}(p)$ so $p \le_{\bbQ[t]} q \Rightarrow
\delta_{\gt}(p) \le \delta_{\gt}(q)$ and 
\sn
\item[$(\beta)$]   if $\ell \in \{0,1\}$ and
$\bar p = \langle p_i:i < \alpha \rangle$ 
is increasing by $\le_{\bbQ[{\gt}]}$ and $\langle\text{ Dom}(S_{p_i}):i
< \alpha\rangle$ is increasing too, $\alpha < \lambda$ a limit ordinal 
and $\delta = \cup\{\Dom(S_{p_u}):i < \alpha\}$ and $i < \alpha
\Rightarrow 0 = S_{p_{i+1}}(\delta_{p_{i+1}})$ and $\ell \in \{0,1\}$
\then \, for some upper bound $p$ of $\bar p,S^{\gt}_p \rest
\delta = \cup \{S^{\gt}_{p_i}:i < \alpha\}$, but $\delta \in S^{\gt}$
and $S^{\gt}_p(\delta) = \ell$; 
moreover for $\ell=0$ there is a $\le_{\bbQ[{\gt}]}$ minimal such $p$ (a
canonical one suffices) with $\delta(p) = \delta +1$
\sn
\item[$(\gamma)$]  if $p \in \bbQ_{\gt}$ then for some $q \in
\bbQ_{\gt}$ we have $p \le_{\bbQ[\gt]} q$ and $\delta_{\gt}(p) <
\delta_{\gt}(q)$ and $S^{\gt}_q(S_{\gt}(p)) = 0$
\sn
\item[$(\delta)$]  so $\name S_{\gt} = \cup\{S^{\gt}_p:p \in
\name G\}$ is a $\bbQ$-name of a function from
$\lambda$ to $\{0,1\}$ usually a stationary subset of $S_{\gt}$
\end{enumerate}
\sn
\item[$(d)$]   first order formulas $\psi_{\gt}(X,Y,\name B),
\varphi_{\gt}(x,y,Y,\name B)$ with $\name B$ a $\bbQ_{\gt}$-name of a
subset of ${\cH}(\lambda),x,y$ individual variables and
$X,Y$ monadic variable, the formulas are in the vocabulary of 
$({\cH}(\lambda),\in,<^*_\lambda)$ 
\sn
\item[$(e)$]   in any universe $\bold V' \supseteq
\bold V$ but with $({}^{\lambda >}\Ord)^{\bold V'} =
({}^{\lambda >}\Ord)^{\bold V}$ (suffice if for some fat stationary
$S \subseteq \lambda,\bold V'$ is gotten from $\bold V$ by forcing by some
strategically $(\lambda \backslash S)$-complete forcing), 
if $\bold G \in \bold V'$ and $\bold G \subseteq \bbQ_{\gt}$ is 
$(< \lambda)$-directed and $\name S_{\gt}[\bold G] 
:= \cup\{S^{\gt}_p:p \in \bold G\}$ is a 
function with domain $\lambda$, and $B =
\name B_{\gt}[\bold G]$ and $A \subseteq {\cH}(\lambda)$ from $\bold V'$ 
 and lastly the formula $(\forall X \subseteq
{\cH}(\lambda)) \psi_{\gt}(X,A,B)$ is satisfied in $\bold V'$ \then \,: 
\sn
\begin{enumerate}
\item[$(*)$]    the formula $\varphi_{\gt} = 
\varphi_{\gt}(x,y,A,B)$ defines in the structure 

$({\cH}(\lambda),
\in,<^*_\lambda,\bold  G)$ a forcing notion $\Bbb Q^{\frak t}_{\bold G,A}$ such
that
\sn
\item[${{}}$]   $(i) \quad$  it is strategically
$(\lambda \backslash \name S_{\gt}[\bold G])$-complete (see \ref{1t.3} below)
\sn
\item[${{}}$]   $(ii) \quad$ for notation simplicity, any
increasing sequence of length $< \lambda$

\hskip25pt  which has an upper bound has a lub (or at least a canonical one)
\sn
\item[${{}}$]   $(iii) \quad \delta = \delta_{{\gt},A}$
is a function from $\bbQ^{\gt}_{\bold G,A}$ to $\lambda$ such that if

\hskip25pt  $\bar p =
\langle p_\varepsilon:\varepsilon < \zeta\rangle$ is an increasing
sequence

\hskip25pt  and $\cup\{\delta(p_\varepsilon):\varepsilon < \zeta\} \in
\lambda \backslash \name S_t[\bold G]$ \then \, $\bar p$ has a lub.
\end{enumerate}
\end{enumerate}
\end{definition}

\begin{definition}
\label{1t.2}
For a $\lambda$-task template $\gt_1$ and $\cW \subseteq \lambda$ we
let $\gt_2 := \gt \rest \cW$ be defined like $\gt_1$ but:
\mn
\begin{enumerate}
\item[$\bullet$]  replacing $S$ by $S \cap \cW$
\sn
\item[$\bullet$]  $\bbQ_{\gt_2} = \bbQ_{\gt_1} \rest \{p \in
\bbQ_{\gt_2}:S_p \subseteq S \cap \cW\}$
\sn
\item[$\bullet$]  similarly to $\varphi,\psi$ (so such $\gt_2$ may not
exist).
\end{enumerate}
\end{definition}

\begin{definition}
\label{1t.3}
1) For regular uncountable $\lambda$, stationary $S \subseteq \lambda$
and ordinal $\alpha \le \lambda$
we say that a forcing notion $\bbQ$ is $\alpha$-strategically
$S$-complete \when \, for every $p^* \in \bbQ$ in the following game the
completeness player has a winning strategy.  A play lasts up to
$\alpha$ moves, in the $\beta$-th move the completeness player chooses
a condition $p_\beta$ and ordinal $\varepsilon_\beta$ such that
$p_\beta$ is an upper bound of $\{p^*\} \cup \{q_\gamma:\gamma <
\beta\}$ and $\varepsilon_\beta$ is an ordinal in $[\bigcup\limits_{\gamma <
\beta} \zeta_\gamma \cup \alpha,\lambda)$ satisfying $[\beta$ limit
$\Rightarrow \varepsilon_\beta = \bigcup\limits_{\gamma < \beta}
\zeta_\gamma\}]$ and the incompleteness player choose $q_\beta$ such
that $p_\beta \le q_\beta \in \bbQ$ and $\zeta_\beta 
\in (\varepsilon_\beta,\lambda)$.  The
completeness player wins if he always has a legal move \underline{or} 
he does not for $\beta$, so $\beta$ is necessarily a limit ordinal but
$\varepsilon_\beta \notin S$.  

\noindent
1A) In part (1) and (2), if $S = \lambda$ and understood from the
context we may omit $S$.

\noindent
2) Let $``(< \alpha)$-strategically $S$-complete" mean
$\beta$-strategically $S$-complete for every $\beta < \alpha$.  If $S =
\lambda$ and $\lambda$ clear from the content we may omit it.  
If we omit ``$(< \alpha)$" we mean $(< \lambda)$.

\noindent
3) If $S \subseteq  \lambda$ is fat, $\bbQ$ is 
$(<\lambda)$-strategically $S$-complete \then \, forcing with $\bbQ$
adds no new sequence from ${}^{\omega >}\bold V$.
\end{definition}

\begin{definition}
\label{1t.5}
1) We call ${\gs}$ a (specific) $\lambda$-task \when \, a triple of
the form 
$(S_{\gs},\psi_{\gs}(X,Y,B_{\gs}),\varphi_{\gs}(x,y,Y,B_{\gs}))$ such that:
\mn
\begin{enumerate}
\item[$(c)'$]   $S_{\gs} \subseteq \lambda$
\sn
\item[$(d)'$]   $\psi_{\gs},\varphi_{\gs}$ are first order formulas as
  in clause (d) of \ref{1t.3}
\sn
\item[$(e)'$]   like clause (e) of Definition \ref{1t.1}(1) except that 
we omit $\bold G$ and use $S$ instead of 
$\name S_{\gt}[\bold G],B_{\gs}$ instead of
$\name B_{\gt}[\bold G]$ and replace $A \mapsto 
\bbQ^{\gt}_{\bold G_{\gs},A}$ defined by $\varphi_{\gt}$ 
by the function $A \mapsto \bbQ^{\gs}_A$ or $\bbQ_{\gs,A}$ defined by
$\varphi_{gs}$.
\end{enumerate}
\mn
2) We say the $\lambda$-task ${\gs}$ is a case of the 
$\lambda$-task template ${\gt}$, in fact is ${\gs}[{\gt},\bold G]$
 for some $(< \lambda)$-directed $\bold G \subseteq \bbQ^{\gt}$ \when \, we
have: ${\gs} = {\gt}<\bold G>$ which means:
\mn
\begin{enumerate}
\item[$(a)$]   $S_{\gs} = \cup\{S^{\gt}_p:p \in \bold G\}$ is
the characteristic function of $S$ (as a subset of $\lambda$)
\sn
\item[$(b)$]   $\varphi_{\gs} = \varphi_{\gt},\psi_{\gt}
= \psi_{\gt},B({\gs}) = B_{\gs} = \name B_{\gt}[\bold G]$.
\end{enumerate}
\mn
3) We say that the $\lambda$-task ${\gs}$ is satisfied
(in a universe $\bold V$) \when \,:
for every $A \subseteq {\cH}(\lambda)$, the $A$-instance of ${\gs}$ 
is satisfied, or ${\gs}$ is satisfied for $A$ which means
that: 

if $({\cH}(\lambda),\in,<^*_\lambda) \models (\forall X
\subseteq {\cH}(\lambda))\psi_{\gs}(X,A,B_{\gs})$ and
$\varphi_{\gt}(x,y,A,B_s)$ define in $({\cH}(\lambda),
\in,\le^*_\lambda)$ a $(< \lambda)$-strategically 
$(\lambda \backslash S_{\gs})$-complete (see Definition \ref{1t.3}) 
forcing notion which we call $\bbQ^{\gs}_A$ 
and ${\cI}_\alpha \subseteq \bbQ^{\gs}_B$ is a predense subset of
$\bbQ^{\gs}_A$ for every $\alpha < \lambda$ 
\then \, $\bbQ_{\gs,A}$ has a directed subset 
not disjoint to ${\cI}_\alpha$ for every $\alpha < \lambda$ (see more in
$(d)^+(\delta)$).
\end{definition}

\begin{remark}
\label{1t.7}
1) Concerning Definition \ref{1t.5} we 
really are interested in finer versions of
being satisfied.  To begin with, $S_{\gs}$ being stationary, and
then preserving further non-existence.  Maybe add to $\lambda$-task
templates (and their cases) a first order formula $\vartheta(X,Y,B)$
and define strong satisfaction (in part (3)) to mean 
that $({\cH}(\lambda),\in) 
\models (\forall X \subseteq {\cH}(\lambda)
[\vartheta(X,\bold G,B)]$.  The idea is that in the forcing
below, for any case of $\gs$ of $\lambda$-task templates, we do not
force more than all the $\bbQ_{\gs,A}$'s.

\noindent
2) So we actually will investigate universes which we get by forcing
as in Definition \ref{1t.9}, \ref{1t.13}(2) below where we iterate the 
forcing from \ref{1t.1}, \ref{1t.5} on all regular cardinals.
\end{remark}

\centerline {$* \qquad * \qquad *$}

\begin{definition}
\label{1t.9}
Assume $\lambda = \lambda^{< \lambda} > \aleph_0$ and $2^\lambda = \lambda^+$.

Let ${\gK}^1_\lambda$ be the class of ${\gp}$ which consists of the
following (so $\bbP_{\gp} = \bbP_{\lambda^+},\bbP_{\gp,\alpha} =
\bbP_\alpha,\bbQ_{\gp,\alpha} = \name{\bbQ}_\alpha$, etc. [similarly
$\gK^1_{\lambda,\xi}$ for $\xi \le \lambda^+$ when we replace
$\lambda^+$ by $\xi$ and omit clause (C)(c) and also $\gK^1_{\lambda,<
  \lambda^+}$])
\mn
\begin{enumerate}
\item[$(A)$]   $(a) \quad \langle \bbP_\alpha,\name\bbQ_\alpha:\alpha
< \lambda^+\rangle$ an iterated forcing with $(< \lambda^+)$ support
(i.e. full)
\sn
\item[${{}}$]   $(b) \quad \bbP = \bbP_{\lambda^+} = 
\cup\{\bbP_\alpha:\alpha < \lambda^+\}$
\sn
\item[${{}}$]   $(c) \quad \bbP_\alpha$ is 
$(< \lambda)$-strategically complete for $\alpha \le \lambda^+$
\sn
\item[$(B)$]   $(a) \quad \name{\bar\cW} = \langle \name W_\alpha:
\alpha < \lambda^+\rangle$
\sn
\item[${{}}$]   $(b) \quad \name\cW_\alpha \subseteq \lambda$ is
increasing modulo ${\cD}_\lambda$ 
\sn
\item[${{}}$]  $(c) \quad$ for limit $\delta,
\name\cW_\delta/\cD_\lambda$ is lub of $\{\name
W_\alpha/\cD_\lambda:\alpha < \delta\}$
\sn
\item[${{}}$]  $(d) \quad f:\lambda^+ \rightarrow \lambda^+$
satisfies $f(\alpha) \le \alpha$ and $f(f(\alpha)) = f(\alpha)$ for
every $\alpha < \lambda^+$
\sn
\item[${{}}$]   $(e) \quad \Vdash_{\bbP_\alpha} ``\bbQ_\alpha$ is 
$(< \lambda)$-strategically $(\lambda \backslash 
\name S_{f(\alpha)})$-complete", see below
\sn
\item[$(C)$]   $(a) \quad$ if $f(\alpha)= \alpha$ then $\bbQ_\alpha = 
\bbQ[\name{\gt}_\alpha]$ so in $\bold V^{\bbP_\alpha}$
 where $\name{\gt}_\alpha$ is a $\bbP_\alpha$-name of

\hskip25pt  a $\lambda$-task template and $\name S_\alpha$ is $\name
S_{\gt_\alpha}$, so a $\bbP_{\alpha +1}$-name of a 

\hskip25pt  subset of $\lambda$ disjoint to $\name\cW_\alpha$ and
 let $\name{\gs}_\alpha =
\gs[\gt_\alpha,\name{\bold G}_{\bbQ_\alpha}]$, is a case of
$\name\gt_\alpha$ 

\hskip25pt  and $\name W_{\alpha +1} = \name W_\alpha \cup
S_{\name{\gs}_\alpha}$ 
\sn
\item[${{}}$]   $(b) \quad$ if $f(\alpha) < \alpha$ then
$\name{\bbQ}_\alpha$ is $\name{\bbQ}_{\name{\gs}_{f(\alpha)},\name
  A_\alpha},\name A_\alpha$ a $\bbP_\alpha$-name of a subset of

\hskip25pt $\cH(\lambda)[G_{\bbP_\alpha}]$, so a $\name A_\alpha$-instance of
$\name{\gs}_\alpha$ which is a case of $\name{\gs}_{f(\alpha)}$
\sn
\item[${{}}$]  $(c) \quad$ if $f(\alpha) = \alpha \le \beta <
\lambda$ and $\name A$ is a $\bbP_\beta$-name of a subset of $H(\lambda)$
\then \, for 

\hskip25pt  unboundedly many $\gamma \in (\beta,\lambda^+)$
if possible, $\name{\bbQ}_\gamma = \bbQ_{\name{\gs}_\alpha,\name A}$.
\end{enumerate}
\end{definition}

\begin{remark}
\label{1t.10}
1)  Why $\name\cW_\alpha$?   Because for $\alpha < \beta \le
\lambda^+,\bbP_\beta / \bbP_\alpha$ is not necessarily $(<
\lambda)$-strategically complete but it is $(< \lambda)$-strategically
$(\lambda \backslash \cW_\alpha)$-complete and $\lambda \backslash
\cW_\alpha$ is fat.

\noindent
2) We may consider task $\gt$ being dealt with in $\beta = f(\beta)$,
   but $S^{\gt}$ being a subset of $S_{\gp,\alpha}$ for some $\alpha =
   f(\alpha)$.  This is quite reasonable but:
\mn
\begin{enumerate}
\item[$(a)$]  to see if $\gt$ is satisfied non-trivially say
  $\Vdash_{\bbP_{\gp}} ``S_{\gp,\beta}$ is stationary" depends also on
  $\gt_\alpha$
\sn
\item[$(b)$]  as for the satisfaction of $\gt_\alpha$, it may be
  subtly changed as on $S_{\gp,\beta}$ we are doing something more.
\end{enumerate}
\mn
3) We may even iterate (2), having a tree structure on $\{\alpha <
\lambda^+:\alpha = f(\alpha)\}$, see \ref{a50}, \ref{a52}, but we
leave it for the time being.
\end{remark}

\begin{definition}
\label{1t.11}
For ${\gp}  \in \gK^1_\lambda$, i.e.  as in \ref{1t.9}.

\noindent
1) We say ${\gp}$ is full \when \,: if $\name{\gt}$ is a 
$\bbP_{\lambda^+}$-name of a $\lambda$-task template with 
$S_{\name{\gt}} \backslash \bigwedge\limits_{\varepsilon < \lambda}
\name S^{\gp}_{\alpha(\varepsilon)}$ stationary for every sequence
$\langle \alpha(\varepsilon):\varepsilon < \lambda\rangle$ of members
of Rang$(f)$, \then \, for unboundedly many $\alpha < \lambda^+,
\gt_\alpha = \name{\gt} \rest (\lambda \backslash \cW_\alpha)$.

\noindent
2) We say ${\gp}$ is $\lambda$-strategically complete \when \, in
 clauses (A)(c),(B)(e) we replace $(< \lambda)$-strategically by
$\lambda$-strategically. 

\noindent
3) Assume $\gT$ is a definition of a set of $\lambda$-task templates.
We say $\gp$ is $\name{\gT}$-full when: 
\mn
\begin{enumerate}
\item[$\bullet_1$]  $\gp \in \gK^1_\lambda$
\sn
\item[$\bullet_2$]  if $f_{\gp}(\alpha) = \alpha$ then
  $\Vdash_{\bbP_{\gp}} ``\name\gt_\alpha[\gp] \in \name\gT"$ 
\sn
\item[$\bullet_3$]  if $\name{\gt}$ is a $\bbP_\alpha$-name and $\alpha
< \lambda^+$ then for unboundedly many $\beta \in (\alpha,\lambda^+)$
we have $\Vdash_{\bbP_\beta}$ ``if $\name\gt \in \gT$ is a
$\lambda$-task templete then $\gt_\beta[\gp] = \name{\gt}"$.
\end{enumerate}
\end{definition}

\begin{definition}
\label{1t.12}
1) We can above use $\gT$, a set of $\le 2^\lambda$ functions, $\gf$ with
domain $\subseteq \gK^1_{\lambda,< \lambda^+}$ such that for $\gp \in
\gK^1_{\lambda,<\lambda^+},\gf(\gp)$ is a $\bbP_{\gp}$-name of a
$\lambda$-task template if defined.

\noindent
2) So $\gp$ is $\gT$-full \when \,:
\mn
\begin{enumerate}
\item[$\bullet_1$]  $\gp \in \gK^1_\lambda$
\sn
\item[$\bullet_2$]  if $f_{\gp}(\alpha) = \alpha$ then 
$\gt_\alpha[\gp] = \gf(\alpha,\bbP,\bar{\bbQ} \rest \alpha)$ 
\sn
\item[$\bullet_3$]  if $\gf \in \gT$ and $\alpha < \lambda^+$ then
either unboundedly many $\beta \in (\alpha,\lambda^+),\gt_\beta[\gp] =
\gf(\gp \rest \beta)$ so the latter is well defined \underline{or}
$\gp \rest \beta \notin \Dom(\gf)$ for every $\beta <
\lambda^+$ large enough.
\end{enumerate}
\end{definition}

\begin{claim}
\label{1t.13}
1) If $\lambda = \lambda^{< \lambda} > \aleph_0$ and $2^\lambda = 
\lambda^+$ \then \, there is a full ${\gp} \in {\gK}_\lambda$.

\noindent
2) Assume in addition that $\gT$ is as in \ref{1t.11}(3) or as in
   \ref{1t.12}, \then \, there is a $\gT$-full $\gp \in \gK^1_\lambda$.

\noindent
3) If $\gp \in \gK^1_{\lambda,\xi}$ \then \, the subset $\bbP'_{\gp}$
   is a dense subset of $\bbP_{\gp}:p \in \bbP_{\gp}$ \underline{iff}:
   for some limit $\delta = \delta_{\gp}(p)$ we have:
\mn
\begin{enumerate}
\item[$(a)$]   $p \in \bbP_{\gp}$
\sn
\item[$(b)$]  if $\alpha \in \Dom(p)$ then $f_{\gp}(\alpha) \in \Dom(p)$
\sn
\item[$(c)$]  if $\alpha = f_{\gp}(\alpha) \in \Dom(p)$ \then
\, $p \rest \alpha$ forces $(\Vdash_{\bbP_{\gp,\alpha}})$ that
$\Dom(p(\alpha)) = \delta \notin S_{\gt_\alpha}$ or $\Dom(p(\alpha)) =
\delta + 1 \wedge p(\alpha)(\delta) = 0 \wedge \delta \in
S_{\gt_\alpha}$ (O.K. even if $\gt_\alpha$ is a $\bbP_\alpha$-name)
\sn
\item[$(d)$]  if $\alpha \in \Dom(p)$ and $f_{\gp}(\alpha) <
\alpha$ \then \, $p \rest \alpha \Vdash_{\bbP_\alpha}
``\delta_{\gt_{f(\alpha)},\name A_\alpha} (p(\alpha)) = \delta$.
\end{enumerate}
\end{claim}

\centerline {$* \qquad * \qquad *$}

We turn to another version of \ref{1t.11}(4)
\begin{claim}
\label{1t.14}
Let $\xi \le \lambda^+$ and $\gp \in \gK^1_{\lambda,\xi}$, see
\ref{1t.9} and $\bbP = \bbP_{\gp}$, etc.

\noindent
1) $\bbP_{\gp}$ is $(< \lambda)$-strategically closed.

\noindent
2) If $\alpha < \xi$ \then \, $\bbP_{\gp}/\bbP_\alpha$ is $(<
   \lambda)$-strategically $\cW^{\gp}_\alpha$-closed.

\noindent
3) If $\gp \in \gK^1_\lambda$ and $f_{\gp}(\alpha) = \alpha$
 \then \, in $\bold V[\bold G]$, every case of the
 $\lambda$-task $\gs[\name{\gt}_\alpha[\gp],
\bold G_{\name{\bbQ}_{\gp,\alpha}}]$ is satisfied.

\noindent
4) In part (3) if $\Vdash_{\bbP_{\gp,\alpha}} ``S^{\gt_\alpha}
   \subseteq \lambda$ is stationary" and $\in \check I[\lambda]$,
   which usually follows then $\Vdash_{\bbP_{\gp,\alpha +1}} ``\name S_\alpha
\subseteq \lambda$ is stationary.

\noindent
5) We can find $\gp \in \gK^1_{\lambda,1}$ such that 
$\Vdash_{\bbP_{\gp,1}} ``\lambda \in \check I[\lambda]"$, hence no
real loss assuming $(\forall \lambda)(\lambda =\cf(\lambda) > \aleph_0
\Rightarrow \lambda \in \check I[\lambda]$.
\end{claim}

\begin{PROOF}{\ref{1t.14}}
Straightforward.
\end{PROOF}

\begin{claim}
\label{1t.15}
Assume $\bold V$ satisfies G.C.H. and $\gT$ is as in \ref{1t.11}(3) or
as in \ref{1t.12}.
\Then \, there is a class forcing notion $\bbP$ such that:
\mn
\begin{enumerate}
\item[$(a)$]   forcing with $\bbP$ preserves cardinality,
cofinality and G.C.H.
\sn
\item[$(b)$]  $\bbP$ is $\bigcup\limits_\lambda \bbP_\lambda$, where
$\langle \bbP_\lambda,\bbQ_\lambda:\lambda$ regular
uncountable$\rangle$, is an iteration, with set support or Easton
support
\sn
\item[$(c)$]   in $\bold V^{\bbP_\lambda}$, the forcing notion
$\bbQ_\lambda$ is as in \ref{1t.13}(1), or every as in \ref{1t.13}(2)
when we have $\gT$
\sn
\item[$(d)$]   for each regular uncountable $\lambda$ the set
${\cH}(\lambda)^{\bold V[\bbP]}$ is equal to ${\cH}
(\lambda)^{\bold V[\bbP_{\lambda}]}$
\sn
\item[$(e)$]  if $\lambda$ is regular uncountable (in $\bold V$),
$\bold V^{\bbP_\lambda} \models ``S \subseteq \lambda$ is stationary" ($\in
\check I[\lambda]$ for simplicity) then $\bold V^{\bbP} \models ``S$ is
stationary and $\diamondsuit_\lambda"$.
\end{enumerate}
\end{claim}

\begin{PROOF}{\ref{1t.15}}
Straight noting that $\bbP/\bbP_\lambda$ is $(<
\lambda)$-st
rategically closed for any $\lambda$ because of
 \ref{1t.14},i.e. of clause (A)(c) of Definition \ref{1t.11}.
\end{PROOF}

\begin{discussion}
\label{1t.19}
1) Claims \ref{1t.13}, \ref{1t.15} may
seem too easy; particularly, if we compare them to more specific cases
from \cite{Sh:587}.  The reason is that the relevant $\name S$'s are
names, so the ``too ambitious" $\lambda$-tasks template
${\gt}$ will not cause the collapse of cardinals but just, e.g. having 
$\name S_\alpha[\gp]$ non-stationary.  Also the various tasks have little
interaction except
that forcing for one $\lambda$-task, create more instances (of $A$'s)
for which we have to force for another $\lambda$-task, and even create
new $\lambda$-tasks templates.
\end{discussion}
\newpage

\section {An example: Relatives of diamonds} \label{anex}

We give an example of a $\lambda$-task template, see \ref{1t.33} on background.

\begin{definition}
\label{1t.31}
Let $\lambda$ be a regular uncountable and $S \subseteq \lambda$ be 
stationary and $\bold c_1,\bold c_2$ are functions from $S$ to the set 
of cardinals $\le \lambda^{++}$ such that $0 < \bold c_1(\delta) <
\bold c_2(\delta)$; if $\bold c_\ell$ is constantly $\kappa_\ell$ we
may write $\kappa_\ell$ instead of $\bold c_i$; if $\delta \in S
\Rightarrow \bold c_2(\delta) = (\bold c_1(\delta))^+$ then we may
omit $\bold c_2$.

We define a $\lambda$-task template 
${\gt} = {\gt}_1(\lambda,S,\bold c_1,\bold c_2)$ by:
\mn
\begin{enumerate}
\item[$\boxplus$]   $S^{\gt} = S$ and $\bbQ_{\gt}$ is defined by:
\begin{enumerate}
\item[$(A)$]  $p \in \bbQ_{\gt}$ iff
\sn
\item[${{}}$]  $(a) \quad p = (\alpha,f,\bar{\cP}) = 
(\alpha_p,f_p,\bar{\cP}_p)$
\sn
\item[${{}}$]  $(b) \quad \alpha < \lambda$
\sn
\item[${{}}$]  $(c) \quad f:\alpha \rightarrow \{0,1\} \text{ such that } 
S^{\gt}_p := f^{-1}(\{1\}) \subseteq S$
\sn
\item[${{}}$]  $(d) \quad \bar{\cP} = \langle {\cP}_\delta:
\delta \in S^{\gt}_p\rangle$
\sn
\item[${{}}$]  $(e) \quad {\cP}_\delta = \cP_{p,\delta} 
\subseteq {}^\delta \delta$
\sn
\item[${{}}$]  $(f) \quad \bold c_1(\delta) \le 
|{\cP}_\delta| < \bold c_2(\delta) \text{ if } \delta \in S^{\gt}_p$,
\sn
\item[${{}}$]  $(g) \quad {\cP}_\delta = \langle f_{\delta,i}:i < \bold
c_p(\delta)\rangle \text{ with no repetitions}$ 
\sn
\item[$(B)$]   order: natural
\sn
\item[$(C)$]  $(a) \quad \name S_{\gt} = \{\alpha <
\lambda:f_p(\alpha)=1$ for some $p \in \name{\bold G}_{\bbQ_{\gt}}\}$,
\sn
\item[${{}}$]  $(b) \quad \name B$ is $\langle \name{\cP}_\delta:\delta \in 
\name S_{\gt}\rangle$, i.e. the set $\{(\delta,{\cP}_\delta)$: for
some $p \in \name{\bold G}_{\bbQ_t}$ we have

\hskip25pt  $f_p(\delta) = 1$ and $\cP_\delta = \cP_{p,\delta}\}$
\sn  
\item[${{}}$]  $(c) \quad \psi_{\gt}(X,Y,\name B) =
\psi_{\gt}(Y,B)$ says that $Y \in 
\prod\limits_{\delta \in \name S_{\gt}} {\cP}_\delta$
\sn   
\item[${{}}$]  $(d) \quad \varphi(x,y,A,B)$ says that for some
ordinal $\alpha < \lambda$:
\sn
\item[${{}}$]   $\qquad (\alpha) \quad x$ has the form $(c_x,f_x)$ with
$c_x \in {}^\alpha 2$ such that $c^{-1}_x\{1\}$ is a 

\hskip40pt closed subset of $\alpha$ and $f_x \in {}^\alpha \lambda$ and
 $\delta \in \name S_{\frak t} \cap c^{-1}_x(\{1\})$

\hskip40pt  $\Rightarrow f_x \restriction \delta 
\notin \name{\cP}_\delta \backslash \{Y(\delta)\}$
\sn  
\item[${{}}$]  $\qquad (\beta) \quad$ similarly $y = (c_y,f_y)$
\sn
\item[${{}}$]  $\qquad (\gamma) \quad c_x \subseteq 
c_y \wedge f_x \subseteq f_y$.
\end{enumerate}
\end{enumerate}
\end{definition}

\begin{remark}
\label{1t.33}
1) So the intention is that:
\mn
\begin{enumerate}  
\item [$(i)$]   $\name S_{\gt}$ is a stationary subset of $S^{\gt}$
\sn
\item[$(ii)$]   $\langle \name{\cP}_\delta:\delta \in 
\name S_{\gt}\rangle$ is a diamond sequence, i.e. $(\forall f \in
{}^\lambda \lambda)(\exists^{\stat} \delta \in \name S_t)(f
\rest \delta \in \name{\cP}_\delta)$ where $\bold c_1(\delta) \le 
|\name{\cP}_\delta| < \bold c_2(\delta)$, alternatively 
${\cP}_\delta = \{\name f_{\delta,i}:i < \bold c(\delta)\}$ where
$\delta \in \name S_t \Rightarrow \bold c_1(\delta) \le \bold
c(\delta) < \bold c_2(\delta)$
\sn
\item[$(iii)$]   if we omit from each $\name{\cP}_\delta$ one 
function for each $\delta \in \name S_{\gt}$ \then \, $(ii)$ stops to hold.
\end{enumerate}
\mn
2) Fleissner proved (Fleisner diamond): assuming
$\bold V = \bold L$ we have: if $\langle {\cP}_\delta:\delta 
\in S\rangle$ satisfies (ii) then for some $\langle
f_\delta:\delta \in S\rangle \in \prod\limits_{\delta \in S}
\name{\cP}_\delta$ the sequence $\langle
f_\delta:\delta \in S\rangle$ is a diamond sequence.

\noindent
3) In \cite{Sh:122} we prove that the assumption 
$\bold V = \bold L$ is necessary and more that is consistently we have
a counterexample (with $|\cP_\delta| > 1$)
\end{remark}

\begin{claim}
\label{1t.35}
Let $\lambda$ be regular uncountable $S \subseteq \lambda$ 
stationary set $\in \check I^{\bold V}[\lambda]$ and $\bold c_1,\bold c_2$ 
are as in Definition \ref{1t.31}.

\noindent
1) \Then \, ${\gt} = {\gt}_1(\lambda,S,\bold c_1,\bold c_2)$ is a
   $\lambda$-task template, see Definition \ref{1t.1}, this holds also
   for $\gt \rest \cW$ when $\cW \subseteq S$ is stationary.
   satisfies the demands \ref{1t.31}.

\noindent
2) If $S \subseteq S^\lambda_\kappa$ is stationary for some 
$\kappa = \text{\rm cf}(\kappa) < \lambda$ is $\bold c_1$ as in
Definition \ref{1t.31} and $\gt = \gt_1(\lambda,S,\bold c)$ and $\bold
c_1(\delta) \le |\delta|$ all this in $\bold V = 
\bold V^{\bbP_{\gq}}$ where $\gq \in \gK^1_\lambda$ is from
\ref{1t.13}(1) (or just $\gT$-full, $\gt_1(\lambda,S,\bold c_1) 
\in \gT$) \then \, for some directed $\bold G \subseteq \bbQ_{\gt}$:
\mn
\begin{enumerate}
\item[$(a)$]   $S_* = \name S_{\gt}[\bold G]$ is a stationary 
subset of $\lambda$ 
\sn
\item[$(b)$]  $\langle \name{\cP}_\delta[\bold G]:\delta
\in \name S_{\gt}[\bold G]\rangle$ is a 
diamond sequence, see \ref{1t.33}(1)(ii)
\sn
\item[$(c)$]  $\cP^*_\delta := \name{\cP}_\delta[\bold G]$
a family of subsets of $\delta$ of cardinality $\ge \bold c_1(\delta)$ but
$< \bold c_2(\delta)$
\sn
\item[$(d)$]   if $Y \in \prod\limits_{\delta \in S_*} \cP^*_\delta$
then for some $f \subseteq {}^\lambda \lambda$ and club $E$ of
$\lambda$ we have $\delta \in S_* \Rightarrow f \rest \delta \notin
\cP^*_\delta \backslash \{Y(\delta)\}$.
\end{enumerate}
\end{claim}

\begin{PROOF}{\ref{1t.35}}
1) Easy.

\noindent
2) So ${\gq} \in \gK^1_\lambda$, in particular $\bar{\bbQ} =
\langle {\bbP}_\alpha,\name{\bbQ}_\beta:\alpha \le \lambda^+,
\beta < \lambda^+\rangle$ are as in \ref{1t.9} and we work in $\bold V$.

So for some $\alpha(0),{\gt}$ is a ${\bbP}_{\alpha(0)}$-name.
Hence by part (1) and \ref{1t.13} for some $\alpha(1) \in
(\alpha(0),\lambda^+)$, we have $f_{\gq}(\alpha) = \alpha$ and
$\name{\gt}^{\gq}_\alpha = {\gt}$ and we choose $\bold G$ the generic
for $\bbQ^{\gq}_{\alpha(1)}$.

Let $\name S_\alpha,\name{\cP}_\delta,\name{\cI}_\delta$ be from the
generic of ${\bbQ}^{\gq}_{\alpha(1)}$.  The problem is to prove that:
\mn
\begin{enumerate}
\item[$(*)$]  $\Vdash ``\langle {\cP}_\delta:\delta \in \name S\rangle$ is a
diamond sequence".
\end{enumerate}
\mn
The proof is as in \cite{Sh:122} or \cite{Sh:587}, but we give
details.

Toward contradiction assume:
\mn
\begin{enumerate}
\item[$(*)_1$]  $p_* \Vdash_{\bbP_{\gp}} ``g \in {}^\lambda \lambda$
  is a counterexample".
\end{enumerate}
\mn
Without loss of generality
\mn
\begin{enumerate}
\item[$(*)$]  $p_* \Vdash \name g \notin ({}^\lambda \lambda)^{\bold
V[\bbP_{\alpha(1)+1}]}$.
\end{enumerate}
\mn
[Why?  Otherwise \wilog \, $p_* \Vdash ``\name g \in
({}^\lambda\lambda)^{\bold V[\bbP_{\alpha(1),2}]}$ and clearly
$\bbQ^{\bbP}_{\alpha(1)}$ naturally forces a $\diamondsuit^*$-sequence.]

Moreover for some $\alpha(3)$ and $p_{**}$ we have 
\mn
\begin{enumerate}
\item[$(*)$]  $(a) \quad \alpha(1) < \alpha(2) \le \alpha(3) < \lambda^+$
\sn
\item[${{}}$]  $(b) \quad p_* \le p_{**} \in \bbP_{\gp}$
\sn
\item[${{}}$]  $(c)(\alpha) \quad$ if $\alpha(2) \ne \alpha(3)$ then
$f_{\gp}(\alpha(2)) = \alpha(1)$ and $p_{**} \Vdash_{\bbP_{\alpha(3)}}
``\name g = \name f_{\gp,\alpha(2)}"$
\sn
\item[${{}}$]  $\quad (\beta) \quad$ if $\alpha(2) = \alpha(3)$ \then \,
for every $\beta \in (\alpha(1),\alpha(3))$ we have 

\hskip35pt $p_{**} \Vdash_{\bbP_{\alpha(2)}} 
``\name g \ne \name f_{\gp,\beta}"$.
\end{enumerate}
\mn
Let
\mn
\begin{enumerate}
\item[$\oplus$]  $\langle e_\alpha:\alpha < \lambda\rangle$ witness $S
\in \check I[\lambda]$
\begin{enumerate}
\item[$(a)$]   $e_\alpha \subseteq \alpha$
\sn
\item[$(b)$]  $\beta \in e_\alpha \Rightarrow e_\beta = e_\alpha \cap
\beta$
\sn
\item[$(c)$]  $\otp(c_\alpha) \le \kappa$
\sn
\item[$(d)$]  $\{\delta \in S:\delta > \sup(c_\beta)\}$ is not
stationary
\sn 
\item[$(e)$]  let $\iota(\alpha) = \otp(c_\alpha)$.
\end{enumerate}
\end{enumerate}
\mn
We define the tree $\cT_i$ as $\{<>\} \cup \{\langle j \rangle:j < 1 +
i\}$.  Now we choose $9\bar p_i,\gamma_i)$ by induction on $i
<\lambda$ such that:
\mn
\begin{enumerate}
\item[$\boxplus(A)$]   $(a) \quad \bar p_i = \langle p_{i,\eta}:\eta
\in \cT_{\iota(\alpha)}\rangle$
\sn
\item[${{}}$]  $(b) \quad p_{i,<>} \in \bbP_{\alpha(2)}$
\sn
\item[${{}}$]  $(c) \quad p_{i,<j>} \in \bbP_{\alpha(3)}$ is above
$p_{**}$
\sn
\item[${{}}$]  $(d) \quad p_{i,<j>} \rest \alpha(2) = p_{i,<>}$
\sn
\item[${{}}$]  $(e) \quad \langle p_{\iota,<j>}:\iota \in c_i\rangle$
is increasing (in $\bbP_{\alpha(3)}$)
\sn
\item[${{}}$]  $(f) \quad p_{i,<>},p_{i,<j>} \in \bbP'_{\alpha(3)}$
moreover
\sn
\item[${{}}$]  $(g) \quad \delta_{\gp}(p_{i,<>}) =
\delta_{\gp}(p_{i,<j>})$ is $\gamma_i$
\sn
\item[$(B)$]  $(a) \quad p_{i,<j>}$ forces a vaue to $\name g \rest
\gamma_i$ call it $g_{i,<j>}$
\sn
\item[${{}}$]  $(b) \quad \langle g_{i,<j>}:j < 1 +i\rangle$ is with
no repetitions
\sn
\item[${{}}$]  $(c) \quad$ if $\beta \in \Dom(p_{i,<j>})$
and $f_{\gp}(\beta) = \alpha(1) < \beta$ \then \, $p \rest \beta$
forces a value to 

\hskip25pt $p(\beta)$ and so its $\in {}^{(\gamma_i)} 2$
\sn
\item[${{}}$]  $(d) \quad$ if $\beta_\ell \in \text{
Dom}(p_{i,<j_\ell>}),f_{\gp}(\beta_\ell) = \alpha(1)$ for $\ell=1,2$
then $p_{i,<j_1>}(\beta_1) \ne$

\hskip25pt $p_{i,<j_2>}(\beta_2)$.
\end{enumerate}
\mn
There is no problem to carry the definition and so $E = \{\delta <
\lambda:\delta$ a limit ordinal and $i < \delta \Rightarrow \gamma_i <
\delta\}$ is a club of $\lambda$.

Let $\bold G_{\alpha(1)} \subseteq \bbP_{\alpha(j)}$ be generic over
$\bold V_\mu$ in $\bold V[\bold G_{\alpha(1)}]$ choose $\delta(*) \in
S \cap E$, so $\otp(C_\delta) = \kappa$ and we know $\kappa_1 = \bold
c_1(\delta),\kappa_2 = \bold c_2(\delta)$ and $\kappa_1 \subseteq
|\delta|$.  We choose $r_1 \in \bold G_{\alpha(1)}$ above $\{p_{i,<>}
\rest \alpha(1):i \in C_\delta\}$ forcing the relevant condition.  For
$\varepsilon < \kappa,g^*_\varepsilon \in {}^\delta 2$ is
$\cup\{g_{i,<\varepsilon>}:i \in C_\delta\}$.

Note that $\langle g^*_\varepsilon:\varepsilon < \kappa_1\rangle$ is a
sequence of members of ${}^\delta \delta$ with no repetitions.  If
$\alpha(2) = \alpha(3)$, then for any $\varepsilon <
\kappa,h_\varepsilon \notin \Lambda_\varepsilon :=
\{\cup\{p_{i,<\varepsilon>}(\beta):i \in C_\delta$ and $\beta \in
\Dom(p_{i,<\varepsilon>}$ and $f(\beta) = \alpha(1) < \beta\}$
and we can choose $r_2 \in \bbP_{\alpha(1)+1}$ above $r_1$ and above
$p_i \rest (\alpha(1)+1)$ for $i \in C_\delta$ such that
$r_2(\alpha(1))(\delta)$ is a set $\cP_\delta$ disjoint to
$\Lambda_\varepsilon$ to which $h_\varepsilon$ belongs.  Easily
$\{r_2\} \cup \{p_{i,<\varepsilon>}:i \in C_\delta\}$ has a common
upper bound and we are done.

So assume $\alpha(2) < \alpha(3)$ hence $\bold f_{\gp}(\alpha(2)) =
\alpha(1)$ and we have $h_\varepsilon =
\cup\{p_{i,<\varepsilon>}(\alpha(2)):i \in C_\delta\}$.  Now find $r_2
\in \bbP_{\alpha(1)+1}$ above $r_1$ and $p_i \rest (\alpha(1)+1)$ for
$i \in C_\delta$ such that $r_2(\alpha(1)(\delta)) =
\{h_\varepsilon:\varepsilon < \kappa_1\}$.  Let $r_3 +
\bbP_{\alpha(2)}$ be above $\{r_2\} \cup \{p_{i,<>} \rest \alpha(4):i
\in C_\delta\}$ clearly exist and \wilog \, $r_3$ forces a value to
$\cY_{\alpha(4)}(\delta)$ say $h_{\varepsilon(*)}$ and now $\{r_3\}
\cup \{p_{i,<\varepsilon(*)>}:i \in C_\delta\}$ has a common upper
bound say $r_4$ and it is as required.
\end{PROOF}

\begin{discussion}
\label{1t.37}
1) What occurs if in \ref{1t.35}(2) we waive ``$S \subseteq
   S^\lambda_\kappa$" and $\bold c_1(\delta) \le \delta$?

We can assume instead:
\mn
\begin{enumerate}
\item[$(*)$]    $(a) \quad S \subseteq \lambda$ is stationary
\sn
\item[${{}}$]  $(b) \quad$ we have $(\alpha)$ or $(\beta)$ where
\begin{enumerate}
\item[${{}}$]  $(\alpha) \quad \kappa = \lambda,S$ is a set of
stronger inaccessible cardinals and $\bold c_1(\delta) \le
2^{|\delta|}$ 

\hskip25pt for $\delta \in S$ 
\sn
\item[${{}}$]  $(\beta) \quad S \subseteq S^\lambda_\kappa,S \in
\check I[\lambda],\kappa = \text{ cf}(\kappa) < \kappa$ and $\bold
c_1(\delta) \le |\delta|^\kappa$, moreover

\hskip25pt  there is a tree $\cT$ with
$\kappa$ levels, $< \lambda$ nodes and 

\hskip25pt $\ge \sup\{\bold c_1(\delta):\delta \in S\},\kappa$-branches.
\end{enumerate}
\end{enumerate}
\mn
The proof of \ref{1t.35} works, with some changes.  First, $\kappa =
\lambda$ then $\delta \in S \Rightarrow \otp(C_\delta) =
\delta$.  Second, having chosen $\bar e = \langle e_\alpha:\alpha <
\lambda\rangle$, we also choose $\langle \cT_\alpha:\alpha <
\lambda\rangle$ such that $\cT_\alpha$ is a tree with
$\otp(C_\alpha)+1$ levels, $<\lambda$ nodes for transparency a sub-tree
of ${}^{\otp(C_\delta)\ge}\delta$ such that $\alpha \in e_\beta
\Rightarrow \cT_\alpha = \cT_\beta \cap
{}^{\otp(C_\alpha)\ge}\lambda$ and $\delta \in S \Rightarrow
|\max(\cT_\delta)| \ge \bold c_1(\delta)$.

Clearly possible in both cases (and we can allow $S$ to be a set of
weakly inaccessible cardinals (so $(\exists \mu < \lambda)(2^\mu =
\lambda))$ and if $\kappa < \lambda$ waive the existence of $\cT$ when
there is such $\bar T$).

Now $\bar p_\alpha$ is $\langle p_{\alpha,\eta}:\eta \in 
\max(\cT_\alpha)\rangle$ and $\beta \in C_\alpha \wedge \eta \in
\max(\cT_\alpha) \Rightarrow p_{\beta,\eta \rest
\otp(C_\beta)} \le p_{\alpha,\eta}$.  In the end having chosen
$\delta \in S \cap E$ we choose a sequence $\langle
\eta_\varepsilon:\varepsilon < \bold c_1(\delta)\rangle$ of pairwise
distinct members of $\max(\cT_\delta)$ and continue as there replacing
$p_{i,<\varepsilon>}$ by $p_{i,\eta_\varepsilon \rest \otp(c_i)}$.
\end{discussion}

\begin{discussion}
\label{1t.39}
1) Is it true that in \ref{1t.35}(2), in $\bold V[\bold G],\bold G
   \subseteq \bbP_{\bold q}$ generic over $\bold V$, we get that for
   every stationary $S \subseteq \lambda$ there are a stationary $S_*
   \subseteq S$ such that the conclusion there holds?  This is almost
   true as if $\alpha(1) < \lambda,\gt^*_2 =
   \gt_{\alpha(1)}[\gq][\bold G]$ is well defined but $\ne \gt$ and
   $S_{\alpha(1)} \subseteq S_*$, then it complicates the forcing
   argument; moreover it may be one which ``promises" $\neg
   \diamondsuit_{\name S_{\alpha(1)}}$.

\noindent
2) However, if we use \ref{1t.13}(2) for $\gT$ consisting of $\gt$
   only, the conclusion aove surely holds.

\noindent
3) What if $\gT$ consists of all $\gt$'s of this form?  In this case
   our framework, i.e. Definition \ref{1t.9} demand that $\langle
   \name S_{\gp,\alpha}:\alpha < \lambda^+,f_{\gp}(\alpha) = \alpha
   \rangle$ have pairwise non-stationary intersections.  We can waive
   this here and then seems O.K.

\noindent
4) So why not generally allow this in \ref{1t.9}?  It is reasonable
   but it complicates things considerably and we do not have an urgent need.
\end{discussion}

\begin{discussion}
\label{1t.41}
A variant of \ref{1t.31} - \ref{1t.35} is:
\mn
\begin{enumerate}
\item[$(\alpha)$]  $(*)_1 \quad \bold c_1(\delta) 
\text{ a cardinal } \le \lambda^+$
\sn
\item[${{}}$]      $(*)_2 \quad \bold c_2(\delta) \text{ a family of subsets
of } \bold c_1(\delta) \text{ closed under subsets}$
\sn
\item[$(\beta)$]  in clause (A) of \ref{1t.31}, clauses (d),(f) are replaced by
\begin{enumerate}
\item[$(d)'$]  $\bar{\cP} = \langle (\cP_\delta,\cJ_\delta):\delta \in
S^{\gt}_p\rangle$
\sn
\item[$(f)$]  ${\cP}_\delta$ is non-empty, $\bold J_\delta \subseteq
\cP(\cP_\delta)$ and $(\cP_\delta,\cJ_\delta)$ is isomorphic to some
pair from $\bold c_2(\delta)$
\end{enumerate}
\item[$(\gamma)$]  in \ref{1t.31}(B)(c)$(\beta)$, $Y$ a function with
domain $\name S_{gt}$ such that $Y(\delta) \subseteq
\cP_\delta,Y(\delta) \notin \cJ_\delta$
\sn
\item[$(\delta)$]  in (B)(c), $f_x \notin \cP_\delta \backslash
Y(\delta)$
\end{enumerate}
\mn
\underline{Older version}:
\mn
\begin{enumerate}
\item[$(\epsilon)$]  in the proof of \ref{1t.35}(2) replace $(*)$ by
\sn
\begin{enumerate}
\item[$(*)$]  and in the stronger version if $\langle 
{\cP}^-_\delta:\delta \in \name{\mathscr{S}}\rangle$ ``if $Y \in 
\prod\limits_{\delta \in \name S} (\name{\cP}_\delta \backslash 
{\cJ}_\delta)$ a $\bbP^{\gq}_{\alpha(i)}$-name $\alpha(1) \le
\alpha(2) < \lambda^+$ then $\langle {\cP}'_\delta:\delta \in \name
S\rangle$ is a diamond sequence.
\end{enumerate}
\end{enumerate}
\end{discussion}

\begin{discussion}
\label{1t.43}
For any $\kappa = \cf(\kappa) < \lambda$ then is a $\lambda$-task
guaranteeing $\diamondsuit^*_{S'_1}$ for some stationary $S_1
\subseteq S$ for any stationary subset of $S$ from $\bold
V^{\bbP_\lambda}$.
\end{discussion}
\newpage

\section {Parameters for completeness of forcing} \label{Parameters}

This section is purely set theoretic.  Trees of conditions continue to
 play major roles, see \cite{Sh:587}; here see the proof of
 \ref{1t.35}, but whereas in the proof of \ref{1t.35} we use
 ``degenerate simple" tree $\cT_i = \{<>\} \cup \{\langle \varepsilon
 \rangle:\varepsilon < i\}$ here we use larger trees and extra
 structure on them.
We concentrate on successor
$\lambda$, for inaccessibles we have to phrase it differently, see \ref{2f.5}.

Recall
\begin{definition}
\label{2f.1}
A filter $D$ on a set $A$ is $(\kappa,\theta)$-regular \when \,
 we can find $A_\alpha \in D$ for $\alpha < \kappa$ such that $i < \mu
\Rightarrow |\{\alpha:i \in A_\alpha\}| < \theta$.
\end{definition}

\begin{definition}
\label{2f.3}
1) For cardinals $\kappa \ge \theta = \cf(\theta)$ and $\delta \ge
\kappa$ a limit ordinal, we 
say $(D,{\cT})$ is a $(\delta,\kappa,\theta)$-special pair \when \,:
\mn
\begin{enumerate}
\item[$(a)$]  $D$ is a $(\kappa,\theta)$-regular filter on $\delta$
\sn
\item[$(b)$]   ${\cT} \subseteq {}^{\delta >}\delta$
 and no $\eta \in {\cT}$ is $\triangleleft$-maximal
\sn
\item[$(c)$]   $<> \in {\cT},{\cT}$ is closed under initial segments
\sn
\item[$(d)$]   ${\cT}$ is closed, that is, if $\eta \in 
{}^{\delta >} \delta$ and $(\forall \alpha < \ell g(\eta))
(\eta \restriction (\alpha +1) \in {\cT})$ then $\eta \in {\cT}$
\sn
\item[$(e)$]   for every $\eta \in \lim_\delta({\cT})$ the set $\{\alpha <
\delta:\text{Suc}_{\cT}(\eta \restriction \alpha)$ is not a singleton$\}$
belongs to $D$.
\end{enumerate}
\mn
1A) We say $\cT$ or $(\cT,D)$ has $\bar c$-successor; if
$\bigwedge\limits_{\eta} c_\eta = 0$ we may omit it
\mn
\begin{enumerate}
\item[$(b)'$]  $\cT$ is $\subseteq {}^{\delta >}\delta$
ordered by $\triangleleft$ with $\eta \in \cT \Rightarrow 
\eta \char 94 \langle c_\eta \rangle \in
\cT$ and we call a $c_\eta$ the default value (for $\eta$) and 
$\bar c = \langle c_\eta:\eta \in \cT \rangle$ is called
the default sequence .
\end{enumerate}
\mn
1B) Let ``$\cT$ is a $\delta$-special tree" 
mean that clause (b),(c),(d) of part (1).

\noindent
2) We say that ${\cT} \subseteq {}^{\delta >} \delta$ 
is $\partial$-lean \Iff \, (b),(c), (d) of part (1):
\mn
\begin{enumerate}
\item[$(f)$]   for every $\alpha < \delta$ the set $\{\eta \in {\cT}:\ell
g(\eta) = \alpha$ and $\Suc_{\cT}(\eta)$ is not a singleton$\}$ 
has $< \partial$ members.
\end{enumerate}
\mn
2A) Saying ``$(D,{\cT})$ is a $(\delta,\kappa,\theta)$-special lean
pair" means $(D,{\cT})$ is a $(\delta,\kappa,\theta)$-special pair
and ${\cT}$ is $\theta$-lean.

\noindent
2B) Saying ${\cT}$ is lean means $\theta$-lean when $\theta$ is
clear from the context.

\noindent
3) For a $(\delta,\kappa,\theta)$-special pair $(D,{\cT})$, let 
\mn
\begin{enumerate}
\item[$(\alpha)$]   $\text{Sub}(D,{\cT}) = \{{\cT}':
{\cT}' \subseteq {\cT} \text{ satisfies clauses (b)-(e) of part (1)
  and } \eta \in {\cT}' \wedge |\Suc_{{\cT}'}(\eta)| > 1
\Rightarrow \Suc_{{\cT}'}(\eta) = \Suc_{\cT}(\eta)\}$,
\sn
\item[$(\beta)$]   $\Ln - \Sub_\partial(D,{\cT}) = 
\{{\cT}' \in \Sub(D,{\cT}):{\cT}' \text{ is } \partial$-lean$\}$
\sn
\item[$(\gamma)$]   in clause $(\beta)$ if $\partial$ is missing and
$(\kappa,\theta)$, or $\theta$ is clear from the context we mean
$\partial = \theta$, so we may write $\Ln-\Sub(D,{\cT})$.
\end{enumerate}
\mn
4) For a subtree ${\cT} \subseteq {}^{\delta >}\delta$, let 
$n \ell({\cT}) = \{\eta \in {\cT}:\ell g(\eta)$ is not a limit ordinal$\}$. 

\noindent
5) We say that $(D,{\cT},\bar E)$ is a $(\delta,\kappa,\theta)$-special
triple or $\delta$-special triple \when \, clauses (a)-(e) of part (1) or
clauses (a)-(e) part (1A) hold and
\mn
\begin{enumerate}
\item[$(g)$]   $\bar E = \langle E_\eta:\eta \in {\cT} \rangle$
\sn
\item[$(h)$]   $E_\eta$ is a filter on $\Suc_{\cT}(\eta)$ or just
a non-empty family of non-empty subsets of $\Suc_{\cT}(\eta)$ 
closed under supersets, that is $u \subseteq v \in \Suc_{\cT}(\eta) 
\wedge u \in E_\eta \Rightarrow v \in E_\eta$.
\end{enumerate}
\mn
6) Let $\Sub(D,{\cT},\bar E)$ be the set of ${\cT}'$ such that 
${\cT}' \subseteq {\cT}$ satisfies clauses (b)-(d) of part (1)
and
\mn
\begin{enumerate}
\item[$(i)$]   for every $\eta \in \lim_\delta(\cT')$ the
set $\{\alpha < \mu:\Suc_{\cT'}(\eta \rest \alpha)
\in E_{\eta \restriction \alpha}\}$ belongs to $D$ where $\Suc_{\cT}(\eta)
:= \{\zeta:(\eta \restriction \alpha) {}^\frown \langle \zeta \rangle 
\in {\cT}\}$.
\end{enumerate}
\mn
7) $\Ln-\Sub_\partial(D,{\cT},\bar E) = \{{\cT}' \in 
\Sub(D,{\cT},\bar E):(D,{\cT}')$ is $\partial$-lean$\}$.

\noindent
7A) We omit $\partial$ in part (7) when $\partial = \theta$ and
$\theta$ is clear from the context.

\noindent
8) We may replace $\delta$ by another set. 
If $\kappa = |\delta|$ we may omit $\kappa$ so write $(\delta,\theta)$; 
usually $\delta$ is a cardinal and then we tend to use $\mu$.
\end{definition}

\begin{remark}
\label{2f.5}
1) For simplifying we usually do not deal with the
$(D,{\cT},\bar E)$-version in this section, but may need it later.

\noindent
2) We may replace $\Suc_{\cT}(\eta) \in E_\eta$ by $\Suc_T(\eta) \ne
\emptyset \mod E_\eta$, i.e.use $E^+_\eta$.
\end{remark}

\begin{definition}
\label{2f.7}
1) We say $\bold p = (D,{\cT},\lambda,S,W,\bar \eta,\bar E)$ is a 
$(\mu,\kappa,\theta)$-parameter \when \, (if $E_\eta =
\{\Suc_{\cT}(\eta)\}$ for $\eta \in {\cT}$ then we may omit $\bar E$):
\mn
\begin{enumerate}
\item[$(a)$]   $(D,{\cT},\bar E)$ is a $(\mu,\kappa,\theta)$-special pair
\sn
\item[$(b)$]   $(\alpha) \quad \lambda$ is regular uncountable
\sn
\item[${{}}$]  $(\beta) \quad S,W$ are stationary subsets of $\lambda$
\sn
\item[${{}}$]   $(\gamma) \quad \lambda \ge \mu$ [? $\lambda =\mu$ possible?]
\sn
\item[${{}}$]  $(\delta) \quad S$ is
$\subseteq S^\lambda_{\cf(\mu)}$ when $\lambda \ge \mu^+$
\sn
\item[${{}}$]  $(\varepsilon) \quad D_{\bold p} := 
\bbD_{\lambda,W}$ is a fat normal filter on ${}^{\lambda
  >}\cP(\lambda)$, see Definition \ref{z23}
\sn
\item[$(c)$]  $(\alpha) \quad \bar \eta = \langle \eta_\delta:\delta \in S
\rangle$
\sn
\item[${{}}$]   $(\beta) \quad \eta_\delta$ is an increasing
sequence of ordinals $< \delta$ from $W$ of limit length

\hskip25pt   with limit $\delta$ for $\delta \in S$ 
\sn
\item[${{}}$]   $(\gamma) \quad \ell g(\eta_\delta) = \mu$ if 
$\mu < \lambda$ and $\eta_\delta \in {}^{\delta \ge}\delta$ otherwise
\sn
\item[${{}}$]  $(\delta) \quad \cup\{\eta_\delta(i):i < \gamma\} 
\notin S$ for any limit $\gamma < \ell g(\eta_\delta),\delta \in S$ 
\sn
\item[$(d)$]  if $\chi$ is large enough, $x \in {\cH}(\chi)$ and $\cY
\in \text{ proj}_{\bbD_{\lambda,W}}(\chi)$, see Definition \ref{z23}

\hskip25pt \then \, we can find 
$\bar N = \langle N_i:i \le \ell g(\eta)\rangle$
and $\delta \in S$ such that
\sn
\item[${{}}$]  $(\alpha) \quad N_i$ is increasing continuous
\sn
\item[${{}}$]   $(\beta) \quad N_{i+1} \cap \lambda = \eta_\delta(i)$
\sn
\item[${{}}$]   $(\gamma) \quad N_i = N_{i+1}$ iff $i$ is a limit
ordinal and $\eta(i) = \cup\{\eta(j):j<i\}$
\sn
\item[${{}}$]  $(\delta) \quad$ if $N_i \ne N_{i+1}$ then $\bar N
\restriction (i+1) \in N_{i+1}$ for $i < \mu$
\sn
\item[${{}}$]   $(\varepsilon) \quad x \in N_i 
\prec ({\cH}(\chi),\in,<^*_\chi)$
\sn
\item[${{}}$]   $(\zeta) \quad N_{i+1} \in \cY$ and $[N_{i+1}]^{<
\theta} \subseteq N_{i+1}$ [2010/4/12 was
$[N_i]^{< \theta} \subseteq N_i$ if $i =0$

\hskip25pt  or $i=j+1$ and $N_i \ne N_j$.]
\end{enumerate}
\mn
2) If $\kappa = \mu$ we may omit $\kappa$, i.e. say $\bold p$ is a
$(\mu,\theta)$-parameter, if $\kappa = \mu,\theta =
\aleph_0$ we may just write $\mu$-parameter. 
\end{definition}

\begin{definition}
\label{2f.9}
Assume $\bold p = (D,{\cT},\lambda,S,W,\bar \eta,\bar E)$ is a
$(\mu,\kappa,\theta)$-parameter.

\noindent
1) We say that a forcing notion $\bbQ$ is $\bold p$-complete \when \,:
\mn
\begin{enumerate}
\item[$(a)$]   $\bold p$ is a $(\mu,\kappa,\theta)$-parameter
\sn
\item[$(b)$]  $\bbQ$ adds no new sequences of ordinals of length $<
\lambda$, moreover is $\bbD_{\bold p}$-complete, see
\ref{2f.7}(1)(b)$(\varepsilon)$ and Definition \ref{z23}
\sn
\item[$(c)$]   there is a function $F$ (called a witness for the 
${\bold p}$-completeness of $\bbQ$) such that for some $x$:
\begin{enumerate}
\item[$\circledast$]   for some $\eta \in \lim_\mu({\cT}')$, 
the set $\{p_{\eta \restriction i}:i < \mu$ non-limit$\}$ 
has an upper bound in $\bbQ$ \when \, :
\sn
\item[${{}}$]   $(\square)\,\,(\alpha) \quad \chi$ is large enough
\sn
\item[${{}}$]  $\qquad (\beta) \quad \langle N_i:i \le \mu \rangle$
as in clause (d) of Definition \ref{2f.7} for 

\hskip45pt  $x' = \langle F,{\bold p},\bbQ \rangle$
\sn
\item[${{}}$]  $\qquad (\gamma) \quad {\cT}' \in 
\text{ Sub}(D,{\cT},\bar E)$
\sn
\item[${{}}$]  $\qquad (\delta) \quad \bar p = \langle p_\eta:\eta \in
n \ell({\cT}') \rangle$ 
\sn
\item[${{}}$]   $\qquad (\varepsilon) \quad i < \mu \Rightarrow
\langle p_\eta:\eta \in n \ell({\cT}') \cap {}^{i \ge}
\mu \rangle \in N_{i+1}$
\sn
\item[${{}}$]   $\qquad (\zeta) \quad$ if $\eta \in n \ell({\cT}')$ and
$\eta \triangleleft \nu \in n \ell({\cT}')$ then 
$p_\eta \le_{\bbQ} p_\nu$ (if defined)
\sn
\item[${{}}$]   $\qquad (\eta) \quad$ the following
set belongs to $D$

\hskip45pt  $\{i < \mu$: the sequence
$\langle p_\eta:\eta \in n \ell({\cT}') \cap {}^{i+1} \mu \rangle$
is equal to 

\hskip45pt  $F(\langle p_\eta:\eta \in n \ell(T') 
\cap {}^{i \ge} {\cT} \rangle)\}$. 
\end{enumerate}
\end{enumerate}
\end{definition}

The following is not directly useful, as we shall use iterations as in
\ref{1t.9} but used as a warmup, so there $S= \name S_{\gt}$ and
$B$ includes $\langle \eta_\delta:\delta \in S\rangle$.
\begin{claim}
\label{2f.11}
1) Assume
\mn
\begin{enumerate}
\item[$(a)$]   $\bold p = (D,{\cT},\lambda,S,W,\bar \eta,\bar E)$ is a
$(\mu,\mu,\theta)$-parameter
\sn
\item[$(b)$]   $\mu = \mu^{< \mu}$ and $\alpha < \mu \Rightarrow
|\alpha|^{< \theta} < \mu$ and $\theta$ is regular
\sn
\item[$(c)$]   $S \subseteq \lambda$ is stationary.
\end{enumerate}
\mn
If $\bar{\bbQ}$ is $(\le \lambda)$-support iteration of $\bold p$-complete
strategically $(\lambda \backslash S)$-complete (or just
$\bbD_{\lambda,W}$-complete) forcing notions \then \,
$\Lim(\bar{\bbQ})$ is $\bold p$-complete and strategically 
$(\lambda \backslash S)$-complete. 

\noindent
2) If $\bbQ$ is $\bold p$-complete and $\bold p$ is
$(\mu,\kappa)$-parameter, \then \, $\bbQ$ 
does not add new $\mu$-sequences of ordinals and
preserve the stationarity of $S^{\bold p}$.
\end{claim}

\begin{PROOF}{\ref{2f.11}}
As in \cite{Sh:587}, FILL?
\end{PROOF}

We now look for ``interesting" filters $D$.
\begin{definition}
\label{2f.13}
1) We say that $D$ is a $(\mu,S,\kappa,\theta,\bold f)$-1-special
filter \Iff \,
\mn
\item[$(a)$]  $S$ is a subset of the uncountable cardinal $\mu$ of
  cardinality $\mu$
\sn
\item[$(b)$]  $D$ is a $(\kappa,\theta)$-regular filter on the
cardinal $\mu$ and $\kappa \ge \theta = \cf(\theta)$
\sn
\item[$(c)$]   $D$ is uniform (so every co-bounded subset of $\mu$
  belongs to $D$) and $S \in D$
\sn
\item[$(d)$]   $D$ is $\theta$-complete
\sn
\item[$(e)$]  $\bold f:\mu \rightarrow \mu \backslash \{0\}$
\sn
\item[$(f)$]   there is a witness $\bar{\cF}$ which means
\sn
\begin{enumerate}
\item[$(\alpha)$]   $\bar{\cF} = \langle {\cF}_\alpha:
\alpha \in S \rangle$
\sn
\item[$(\beta)$]   $\cF_\alpha \subseteq \prod\limits_{i < \alpha}
  \bold f(i)$ \underline{or}
\sn
\item[$(\gamma)$]   $|\cF_\alpha| < \theta$
\sn
\item[$(\delta)$]    for every $f \in \prod\limits_{i < \mu}
\bold f(i)$ the set
$\{\alpha < \mu:f \restriction \alpha \in {\cF}_\alpha\}$ belongs to $D$.
\end{enumerate}
\mn
1A) Omitting $\bold f$ we mean $i < \mu \Rightarrow \bold f(i) = \theta$.

\noindent
2) We say $(D,{\cT})$ is $(\mu,S,\kappa,\theta,\bold f)-2$-special 
pair \when \,:
\mn
\begin{enumerate}
\item[$(A)$]   $\mu \ge \kappa > \theta,\theta$ regular, $S \subseteq
\mu = \sup(S)$ and $\bold f:\mu \rightarrow \mu +1$
\sn
\item[$(B)$]   $(a) \quad D$ is a $(\mu,S,\kappa,\theta,\bold
  f)-1$-filter, \underline{older}: $(\kappa,\theta)$-regular uniform 
filter on $\mu$
\sn
\item[${{}}$]   $(b) \quad D$ is $\theta$-complete and $S \in D$
\sn
\item[${{}}$]  $(c) \quad \cT$ is a subtree of $\bigcup\limits_{\alpha
< \delta} \prod\limits_{\beta < \alpha} \bold f(\beta) \subseteq
{}^{\mu >}\mu$ of cardinality $\le \mu$
\sn
\item[${{}}$]  $(d) \quad {\cT}$ is a $\theta$-lean, see
Definition \ref{2f.3}(2) and ${\cT} \subseteq {}^{\mu >} \mu$ 
\sn
\item[${{}}$]    $(e) \quad (D,{\cT})$ is a
$(\mu,\kappa,\theta)$-special pair
\sn
\item[${{}}$]    $(f) \quad$ if $\eta \in {\cT}$ and suc$_{\cT}(\eta)$
is not a singleton then $\{\varepsilon:\eta \char 94 \langle
\varepsilon \rangle \in \cT\} =$

\hskip25pt $\bold f(\ell g(\eta))$ (if we have
$\bar E$ as in \ref{2f.3}(5) then it is natural to demand

\hskip25pt  just that the set to $\in E_\eta$).
\end{enumerate}
\mn
3) We say $(D,\bar u,\bold f,\bar v)$ is $(\mu,S,\kappa,\theta,\bold
f)$-3-spcial \when \, $\cT,\cT_*,f_*$ satisfies: (omitting $\bar v$
means for some $\bar v$, omitting $\bar u,\bar v$ means for some 
$\bar u,\bar v$)
\mn
\begin{enumerate}
\item[$(A)$]  $(a) \quad \mu \ge \kappa > \theta = \cf(\theta)$
\sn
\item[${{}}$]  $(b) \quad S \subseteq \mu$
\sn
\item[$(B)$]  $(a) \quad D$ is a $(\kappa,\theta)$-regular
\sn
\item[${{}}$]  $(b) \quad D$ is $\theta$-complete
\sn
\item[${{}}$]  $(c) \quad \cT \subseteq {}^{\mu >}\mu$ is as in part
(2), hence $\theta$-lean
\sn
\item[${{}}$]  $(d) \quad (D,\cT)$ is a $(\mu,\kappa,\theta)$-special pair
\sn
\item[${{}}$]  $(e) \quad S \in D$
\sn
\item[$(C)$]  $(a) \quad \bar u = \langle u_\alpha:\alpha <
\mu\rangle$ is $\subseteq$-increasing with union $\mu$
\sn
\item[${{}}$]  $(b) \quad \Dom(\bold f_*) = \mu$ and $\bold f_*(\alpha)
= {}^{(u_\alpha)}\bold f(\alpha)$
\sn
\item[${{}}$]  $(c) \quad \bar v = \langle v_\alpha:\alpha <
\mu\rangle,v_\alpha \subseteq \mu,|v_\alpha| < \theta$
\sn
\item[${{}}$]  $(d) \quad \bar{\cF} = \langle
(\cF_\alpha,<_\alpha:\alpha \in S\rangle$
\sn
\item[${{}}$]  $(e) \quad \cF_\alpha \subseteq \{\bar\eta:\bar\eta =
\langle \eta_i:i \in v_\alpha\rangle$ and $\eta_i \in \cT$ and $\ell
g(\eta_i) = \alpha$ for $i \in v_\alpha\}$
\sn
\item[${{}}$]  $(f) \quad |\cF_\alpha| < \theta$
\sn
\item[${{}}$]  $(g) \quad <_\alpha$ is a well ordering of $v_\alpha$
\sn
\item[${{}}$]  $(h) \quad$ if $\eta_i \in \lim_\mu(\cT)$ for $i < \mu$
\then \, $\{\alpha < \mu:\langle \eta_i \rest \alpha:i \in v_\alpha
\rangle \in \cF_\alpha\} \in D$
\sn
\item[${{}}$]  $(i) \quad$ if $<_*$ is a well ordering of $\lambda$ then
$\{\alpha:<_\alpha = <_* \rest v_\alpha\} \ne \emptyset \mod D$
\sn
\item[${{}}$]  $(j) \quad \cT_*$ is defined by: $\eta \in \cT_*$ iff:
\sn
\begin{enumerate}
\item[${{}}$]   $(\alpha) \quad \eta$ is a sequence
\sn
\item[${{}}$]  $(\beta) \quad \ell g(\eta) < \mu$
\sn
\item[${{}}$]  $(\gamma) \quad \eta(\alpha) = \langle
\eta(\alpha,\varepsilon):\varepsilon \in v_\alpha \rangle$ and
$\eta(\alpha,\varepsilon) \in \cT$ has length $\alpha + 1$
\sn
\item[${{}}$]  $(\delta) \quad$ if $\alpha < \beta < \ell g(\eta)$ and
$\varepsilon \in u$ then $\eta(\alpha,\varepsilon) \triangleleft
\eta(\beta,\varepsilon)$
\sn 
\item[${{}}$]  $(\varepsilon) \quad$ if $\alpha < \beta < \ell
g(\eta)$ and $\varepsilon \in u_\beta \backslash u_\alpha$ then
$\eta(\beta,\varepsilon)(\alpha)$ is 0 

\hskip25pt or just the default
\sn
\item[${{}}$]  $(\zeta) \quad$ if $\alpha < \ell g(\eta),\varepsilon
\in u_\alpha \backslash v_\alpha$ \then \,
$\eta(\alpha,\varepsilon)(\alpha)$ is 0 or just the default.
\end{enumerate}
\end{enumerate}
\mn
4) We can above allow $\mu = \aleph_0 = \kappa = \theta$.

\noindent
5) If we omit $\kappa$ we mean $\mu = \kappa$.
\end{definition}

\begin{claim}
\label{2f.14}
1) Assume $D$ is a $(\mu,S,\kappa,\theta,\bold f)-1$-special filter,
$\kappa = \mu$ and $\alpha < \mu \Rightarrow |\prod\limits_{\beta
  <\alpha} \bold f(\beta)| < \mu$.
\Then \, for some $\cT$, the pair $(D,\cT)$ is
$(\mu,S,\kappa,\theta,\bold f)$-2-special.

\noindent
2) Assume $\bold f:\mu \rightarrow \mu +1$ and $\bar u,\bold f_*$ are
   as in \ref{2f.13}(3)(C),(a),(b).  If $D$ is
$(\mu,S,\kappa,\theta,\bold f_*)-1$-special \then \, $D$ is
   $(\mu,S,\kappa,\theta,\bold f)-3$-special.
\end{claim}

\begin{PROOF}{\ref{2f.14}}
1) Let $\bar{\cF}$ be a witness for $D$ being
   $(\mu,S,\kappa,\theta,\bold f)-1$-special filter, i.e. as in
   \ref{2f.13}(2).  We define $\cT$ as the set of $\eta$ such that:
\mn
\begin{enumerate}
\item[$(*)_1$]  $(a) \quad \eta \in \prod\limits_{i < \alpha} \bold
  f(i)$ for some $\alpha < \mu$
\sn
\item[${{}}$]  $(b) \quad$ if $\alpha < \ell g(\eta)$ and $\eta \rest
  \alpha \notin \cF_\alpha$ (which holds trivially if $\alpha \notin
  S$) then $\eta(\alpha)=0$.
\end{enumerate}
\mn
2) Should be clear.
\end{PROOF}

\noindent
For $\mu$ regular below we construct such filters.
\begin{claim}
\label{2f.15}
If $\mu > \theta$ are regular, $S \subseteq \mu$ is stationary 
such that $\diamondsuit_S$ holds and $\bold f:\mu 
\rightarrow \mu$ satisfies $i < \mu \Rightarrow \theta \le \bold f(i)
\le i$ \then \, there some $D$ is a 
$(\mu,S,\mu,\theta,\bold f)-1$-special filter.
\end{claim}

\begin{PROOF}{\ref{2f.15}}
As $\diamondsuit_S$ clearly $\mu =
\mu^{< \mu}$ hence $\mu = \mu^{< \theta}$.
Let cd:${}^{\theta >} \mu \rightarrow
\mu$ be one to one onto and cd$(\langle \alpha_i:i <j \rangle) \ge
\sup\{\alpha_i:i <j\}$ and let $g:\mu \rightarrow \theta$ and
cd$_i:\mu \rightarrow \mu$ be such that $\eta \in {}^{\theta >} \mu
\Rightarrow g(\text{cd}(\eta)) = \ell g(\eta)$ and 
$\bigwedge\limits_{i < \ell g(\eta)} \eta(i) = \text{ cd}_i(\text{cd}(\eta))$.

As $\diamondsuit_S$ holds, there is a sequence $\langle
f_\delta:\delta \in S \rangle$ such that
\mn
\begin{enumerate}
\item[$(i)$]    $f_\delta \in {}^\delta \delta$
\sn
\item[$(ii)$]   if $f \in {}^\mu \mu$ then the set $\{\delta \in
S:f_\delta = f \restriction \delta\}$ is a stationary subset of $\mu$.
\end{enumerate}
\mn
Now for $\delta \in S$ we define ${\cF}_\delta \subseteq [{}^\delta
\delta]^{< \theta}$ as follows:

\begin{equation*}
\begin{array}{clcr}
{\cF}_\delta = \{f \in {}^\delta \delta:&\text{ for some } j<i <
\theta \text{ for every } \alpha < \delta \text{ we have } \\
  &g(f_\delta(\alpha)) = i \text{ and } 
f(\alpha) = \text{ cd}_j(f_\delta(\alpha))\}.
\end{array}
\end{equation*}

\mn
By the choice of $g$ we have 
$g(f_\delta(0)) < \theta$ so clearly ${\cF}_\delta \in [{}^\delta
\delta]^{< \theta}$.  Now if $\langle g(f_\delta(\alpha)):\alpha <
\delta\rangle$ is not constant then ${\cF}_\delta = \emptyset$, so we
can ignore such $\delta$'s.

Now for every $h \in {}^\mu \mu$ let

\[
A_h = \{\delta \in S:h \restriction \delta \in {\cF}_\delta\}
\]

\mn
and let $D$ be the $\theta$-complete filter on $\mu$ generated by
$\{A_h:h \in {}^\mu \mu\} \cup \{C:C \text{ a club of } \mu\}$.  Now
\mn
\begin{enumerate}
\item[$(\alpha)$]   $\emptyset \notin D$.
\end{enumerate}
\mn
[Why?  It suffices to prove that $\cap\{A_{h_\varepsilon}:\varepsilon
< \zeta_1\} \cap \{C_\varepsilon:\varepsilon < \zeta_2\}$ is non-empty
when $\zeta_1,\zeta_2 < \theta,h_\varepsilon \in {}^\mu \mu$ for
$\varepsilon < \zeta_1$ and $C_\varepsilon$ is a club of $\mu$ for
$\varepsilon < \zeta_2$.  Let $f \in {}^\mu \mu$ be defined by $f(\alpha) =
\text{ cd}(\langle h_\varepsilon(\alpha):\varepsilon < \zeta_1
\rangle) \in \mu$, hence $S_f = \{\delta \in S:f \restriction \delta =
f_\delta$ (hence $\in {}^\delta \delta)\}$, is a stationary subset of
$\mu$.  Clearly
\mn
\begin{enumerate}
\item[$\bullet$]   $\alpha < \delta \in S_f \Rightarrow g(f_\delta(\alpha)) =
\zeta_1$ and $\bigwedge\limits_{\varepsilon < \zeta_1} 
\text{ cd}_\varepsilon(f_\delta(\alpha) = h_\varepsilon(\alpha)$ and
\sn
\item[$\bullet$]   $\cap\{C_\varepsilon:\varepsilon < \zeta_2\}$ 
is a club of $\mu$.
\end{enumerate}
\mn
So there is $\delta \in S_f \cap
\bigcap\{C_\varepsilon:\varepsilon < \zeta_2\}$ and clearly
$\varepsilon < \zeta_1 \Rightarrow h_\varepsilon \restriction \delta
\in {\cF}_\delta \Rightarrow \delta \in A_{h_\varepsilon}$ so
$\delta \in \bigcap\limits_{\varepsilon < \zeta_1} A_{h_\varepsilon} \cap
\bigcap\limits_{\varepsilon< \zeta_2} C_\varepsilon$ hence we are done.]
\mn
\begin{enumerate}
\item[$(\beta)$]   $D$ is $(\mu,\theta)$-regular. 
\end{enumerate}
\mn
[Why?  For $\varepsilon < \mu$ let $h_\varepsilon \in {}^\mu \mu$ be
constantly $\varepsilon$; so $\{A_{h_\varepsilon}:\varepsilon < \mu\}
\subseteq D$ and no $\alpha < \mu$ belongs to $\ge \theta$ of them.]

So clearly clauses $(a) - (f)$ of Definition \ref{2f.13} hold.
\end{PROOF}
\bigskip

\noindent
Also we can combine such filters.
\begin{claim}
\label{2f.17}
$D$ is a $(\mu,S,\kappa,\theta,\bold f)-1$-special filter \when \, for
some $S,\bar D,\bar\mu,\bar\kappa,\bar S,\bar\kappa,\bar{\bold f}$ we have:
\mn
\begin{enumerate}
\item[$(a)$]   $\mu$ is uncountable and $\theta = cf(\theta) 
\le \kappa \le \mu$
\sn
\item[$(b)$]   $S_* \subseteq \mu$ is unbounded
\sn
\item[$(c)$]   $\bar \mu = \langle \mu_\delta:\delta \in S_* \rangle$
such that for every $\delta \in S_*$ we have:
\sn
\begin{enumerate}
\item[$(\alpha)$]   $\delta \le \mu_\delta < \mu$
\sn
\item[$(\beta)$]  $\alpha \in S_* \Rightarrow 
\kappa_\alpha \ge \prod\limits_{\delta < \alpha} \,
\prod\limits_{\alpha < \mu_\delta} \bold f_\delta(\alpha)$
\sn
\item[$(\gamma)$]  $(\delta,\delta + \mu_\delta) \cap S_* = \emptyset$
and [nec?] $\delta = \otp(\delta \backslash \cup\{[\delta',\delta' +
\mu_\delta):\delta' \in S \cap \delta\})$  for $\delta \in S_*$]
\end{enumerate}
\sn
\item[$(d)$]   $D_\delta$ is a
$(\mu_\delta,S_\delta,\kappa_\delta,\theta_\delta,\bold
f_\delta)-1$-special filter, see \ref{2f.13}(1) and $\theta_\delta \le
\theta,\kappa_\delta < \kappa$
\sn
\item[$(e)$]   $D_\mu$ is a $\theta$-complete
$(\cf(\mu))$-regular filter on $\mu$ such that $S_* \in D$ containing
  the co-bounded subsets of $\mu$
\sn
\item[$(f)$]   if $\theta_* < \theta$ and $\kappa_* < \kappa$ then
$\{\alpha \in S_*:\theta \ge \theta_\alpha > \theta_*$ and $\kappa_\alpha >
\kappa_*$ and $\kappa_\alpha < kappa\} \in D$
\sn
\item[$(g)$]   we define the function $\bold f$ with domain $\mu$ such
that $\bold f(\delta + \alpha) = \bold f_\delta(\alpha)$ if $\delta
\in S_*,\alpha < \mu_\delta$ and $\bold f$ is 1 otherwise
\sn
\item[$(h)$]    $D = \{A \subseteq \mu$: the set $\{\delta
\in S_*:\{\alpha < \mu_\delta:\delta + \alpha \in A\} \in D_\delta\}$
belongs to $D_\mu\}$
\sn
\item[$(i)$]    $S = \cup\{\delta + \alpha:\delta \in S_*$ and $\alpha
\in S_\delta\}$.
\end{enumerate}
\end{claim}

\begin{PROOF}{\ref{2f.17}}

\noindent
\underline{Clause (a) of \ref{2f.13}}:  By clause (b) + (c) we have
$\mu = \sup\{\mu_\delta:\delta \in S_*\}$.  By clause (d), if $\delta
\in S_*$ then $S_\delta$ is an unbounded $\gamma$ subset of
$\mu_\delta$.  Hence by clauses (b) + (i) clearly $S$ is an unbounded
subset of $\mu$ and even of cardinality $\mu$.

\noindent
\underline{Clause (b) of \ref{2f.13}}:  First, $D$ is a filter on $D$ by
its definition (in clause (h) here) as $D_\mu$ is a filter on $\mu$ to
which $S_*$ belongs (see clause (e) of our assumptions) and
$D_\delta$ is a filter on $\mu_\delta$ for $\delta \in S_*$.  
Second, why is $D$ $(\kappa,\theta)$-regular?  Recall $D_\delta$ is
$(\kappa_\delta,\theta_\delta)$-regular by clause (d) of the
assumption and let $\langle A_{\delta,i}:i < \kappa_\delta\rangle$
witness it.

Let $A_\varepsilon = \{\delta + \alpha:\kappa^{j(\varepsilon)} <
\kappa_\delta$ and $\alpha < \mu_\delta$ and $\alpha \in
A_{\delta,\varepsilon}\}$ so it suffices to show that $\bar A = \langle A_i:i
< \kappa\rangle$ witness $D$ is $(\kappa,\theta)$-regular.

First, if $\varepsilon < \kappa$ then 
\mn
\begin{enumerate}
\item[$\bullet$]  $A^1_\varepsilon = \{\delta \in
  S_*:\kappa^{j(\varepsilon)} < \kappa_\delta\} \in D_\mu$ by clause
  (f) of the assumption
\sn
\item[$\bullet$]  if $\delta \in A^1_\varepsilon$ then $\{\alpha <
  \mu_\delta:\delta + \alpha \in A_\varepsilon\} =
  A_{\delta,\varepsilon} \in D_\delta$
\end{enumerate}
\mn
so together $A_i \in D$.

Next let $v \subseteq \kappa,|v| \ge \theta$ and we should prove that
$\cap\{A_i:i \in v\} = \emptyset$.  Toward contradiction assume
$\delta + \alpha \in \cap\{A_\varepsilon:\varepsilon \in v\}$ where
$\delta \in S_*,\alpha < \mu_\delta$, hence clearly $\varepsilon \in v
\Rightarrow \varepsilon < \kappa_\delta$ 
and $\varepsilon \in v \Rightarrow \alpha \in A_{\delta,\varepsilon}$
hence $\{i <\kappa_\delta:\alpha \in A_{\delta,i}\} \supseteq v$ has
cardinality $\ge\theta$, contradiction to the choice of $\langle
A_{\delta,i}:i < \kappa_\delta\rangle$, so indeed $\cap\{A_i:i \in v\}
= \emptyset$ so $\bar A$ exemplifies $D$ is $(\kappa,\theta)$-regular.

\noindent
\underline{Clause (c) of \ref{2f.13}}:  Why $D$ is a uniform filter on
$\mu$?

Assume $\chi < \mu$ and $A \in D$ and we shall prove $|A| \ge r$.  Let
$A_1:\{\delta \in S_*:A_{2,\delta} \in D_\mu\}$ where $A_{2,\delta} =
\{\alpha < \mu_\delta:\delta + \alpha \in A_i\}$.  Now $D_\mu$
contains the co-bounded subsets of $\mu$, hence $\mu = \sup(A_1)$ so
choose $\delta \in A_1$ such that $\delta > \chi$.  Now
$|A_{2,\delta}| = \mu_\delta$ by \ref{2f.13}(1)(b), the uniformity and
clause (d) of the assumption.  Also $|A \rest [\delta,\mu_\delta +
\delta)| = |A_{2,\delta}|$ and $\mu_\delta \ge \delta$ by clause (c)
  of the assumption.  Together $|A| \ge |A \cap [\delta,\delta +
    \mu_\delta)| = |A_{2,\delta}| \ge \mu_\delta \ge \delta \ge \chi$
    so we are done.

Also $S \in D$ by the choice of $D$ in
clause (h) and the choice of $S$ in clause (i).

\noindent
\underline{Clause (d) of \ref{2f.13}}:  $D$ is $\theta$-complete as
\mn
\begin{enumerate}
\item[$\bullet$]  $D_\mu$ is $\theta$-complete by clause (e) of the
  assumption
\sn
\item[$\bullet$]  $D_\delta$ is $\theta_\delta$-complete
by clause (d) of the assumption
\sn
\item[$\bullet$]  $\theta = \cf(\theta)$ by clause (a)
\sn
\item[$\bullet$]  \underline{by clause (f)},
$\theta_* < \theta \Rightarrow \{\delta \in S_*:\theta_*
< \theta_\delta\} \in D$
\end{enumerate}
\mn
and the choice of $D$ in clause (h).

\noindent
\underline{Clause (e) of \ref{2f.13}}:  By the choice of the function
$\bold f$ in clause (g), indeed $\bold f:\mu \rightarrow \mu \backslash
\{0\}$.

\noindent
\underline{Clause (f) of \ref{2f.13}}:  For each $\delta \in S_*$ let
$\langle g_{\delta,\varepsilon}:\varepsilon < \kappa_\delta\rangle$
list $\prod\limits_{\alpha \in S_* \cap \delta} \, \prod\limits_{i <
\mu_\alpha} \bold f_\alpha(i)$, possible by clause $(c)(\beta)$ of the
assumption.  For each $\delta \in S_*$ by the assumption ``$D_\delta$ is
$(\mu_0,S,\kappa_\delta,\theta_\delta,\bold f_\delta)-1$-special 
filter", there is a sequence $\bar A_\delta =
\langle A_{\delta,\varepsilon}:\varepsilon < \kappa_\delta\rangle$
exemplifying ``$D_\delta$ is $(\kappa_\delta,\theta_\delta)$-regular",
and a sequence $\bar\cF_\delta = \langle \cF_{\delta,\alpha}:\alpha \in
S_\delta\rangle$ witnessing clause (f) of \ref{2f.13}(1) for this tuple and let
$\cF'_{\delta,\alpha} = \{f:\text{ Dom}(f) = [\delta,\delta + \alpha)$
and $\langle f(\delta + \alpha):\alpha < \mu_\delta\rangle \in
\cF_{\delta,\alpha}\}$ is a subset of $\Pi\{\bold f(\beta):\beta \in
[\delta,\delta + \alpha)\}$ of cardinality $< \theta_0$.

We define $\bar\cF = \langle \cF_\beta:\beta \in S\rangle$ by
\mn
\begin{enumerate}
\item[$(*)$]  if $\beta = \delta + \alpha,\delta \in S_*$ and $\alpha
\in S_\delta$ then $\cF_\beta = \{f \cup g:f \in \cF'_{\delta,\alpha}$
and $g \in \{g_{\delta,\varepsilon}:\alpha \in
A_{\delta,\varepsilon}\}\}$.
\end{enumerate}
\mn
We should check the four subclause of clause (f) of \ref{2f.13}.

\noindent
\underline{Subclause $(\alpha)$}: holds; trivially by the choice of
$\bar{\cF}$.

\noindent
\underline{Subclause $(\beta)$}:

Recall the choice of $\bold f$ and $\cF_\alpha$.  

\noindent
\underline{Subclause $(\gamma)$}:  As $\delta \in S \Rightarrow
\theta_\delta \le \theta$ it follows that if $\beta = \delta + \alpha
\in S,\delta \in S_*,\alpha \in S_\delta$ then $|\cF_\beta| =
|\cF_{\delta,\alpha}| \times |\{\varepsilon:\alpha \in
A_{\delta,\varepsilon}\}|$ a product of two cardinals $< \theta_\delta
\le \theta$, as required in subclause $(\gamma)$.

\noindent
\underline{Subclause $(\delta)$}:  Let $f \in \prod\limits_{\alpha <
  \mu} \bold f(\alpha)$, for $\delta \in S_*$ let $f_\delta \in
\prod\limits_{\alpha < \mu_\delta} \bold f_\delta(\alpha)$ be 
defined by $f_\delta(\alpha) =
f(\delta + \alpha)$ for $\alpha \in S_\delta$.  So for every $\delta
\in S$, the set $B_\delta = \{\alpha < \mu_\delta:f_\delta \rest
\alpha \in \cF_{\delta,\alpha}\} = \{\alpha < \mu_\delta:f \rest
[\mu_\delta,\mu_\delta + \alpha) \in \cF'_{\delta,\alpha}\}$ belongs
to $D_\delta$.  Also $f_\delta \rest \delta \in
\{g_{\delta,\varepsilon}:\varepsilon < \kappa_\delta\}$, so we can
choose $\varepsilon(\delta) <\kappa_\delta$ such that
$g_{\delta,\varepsilon} = f_\delta \rest \delta$.  Hence if $\alpha <
\mu_\delta$ belongs to $B_\delta \cap A_{\delta,\varepsilon(\delta)}$
then $f \rest [\delta,\delta + \alpha) \in \cF'_{\delta,\alpha},f \rest
\delta \in\{g_{\delta,\varepsilon}:\alpha \in
A_{\delta,\varepsilon}\}$ together $f \rest (\delta + \alpha) \in
\cF_{\delta + \alpha}$.

By the definition of $D$ it follows that $\{\alpha \in S:f \rest
\alpha \in \cF_\alpha\}$ belongs to $D$ as required.
\end{PROOF}

\begin{conclusion}
\label{2f.20}
Assume $\mu > \sigma,\bold f \in {}^M\{\sigma\}$ and $(\forall \kappa
< \mu)(2^\kappa = \kappa^+)$.

\noindent
1) There is a $(\mu,\mu,\mu,\aleph_0,\bold f)-1$-special filter.

\noindent
2) Moreover, there is $(\mu,\mu,\mu,\aleph_0,\bold f)-2$-special pair $(D,I)$.
\end{conclusion}

\begin{proof}
1) By \ref{2f.15} and \ref{2f.17}.

\noindent
2) By part (1) and \ref{2f.14}(1).
\end{proof}

\newpage

\bibliographystyle{alphacolon}
\bibliography{lista,listb,listx,listf,liste,listz}

\end{document}